\documentclass[12pt]{article}
\usepackage{amsmath, amssymb, a4, tabbert}
\usepackage{graphics}
\newcommand{\car}{{\cal A}_r(\om)}
\newcommand{\hw}{\widehat{W}}
\newcommand{\wrar}{\rightharpoonup}
\newcommand {\supp}{\text{supp}}
\newcommand {\fou}{{\cal F}}
\newcommand {\finv}{{\cal F}^{-1}}
\newcommand{\tb}{\textbf{b}}
\begin{document}
\bibliographystyle{alpha}
\title{Resolvent Estimates and Maximal Regularity in Weighted Lebesgue Spaces of
the Stokes Operator \\in Unbounded Cylinders}
\author{Myong-Hwan Ri and Reinhard Farwig}
\date{}
\maketitle
\begin{abstract}\noindent
We study resolvent estimate and maximal regularity of the Stokes operator in $L^q$-spaces with
exponential weights in the axial directions of unbounded cylinders of $\R^n,n\geq 3$.
For straights cylinders we obtain these results in Lebesgue spaces with exponential weights in the axial direction
and Muckenhoupt weights in the cross-section.
Next, for general cylinders with several exits to infinity
we prove that the Stokes
operator in $L^q$-spaces with exponential weight along the axial directions
generates an exponentially decaying analytic semigroup and has maximal regularity.

The proofs for straight cylinders use an operator-valued Fourier multiplier theorem and
techniques of unconditional Schauder decompositions based on the
${\mathcal R}$-boundedness of the family of solution operators
for a system in the cross-section of the cylinder parametrized by the phase variable of the one-dimensional
partial Fourier transform.
For general cylinders we use
cut-off techniques based on the result for straight cylinders and
the result for the case without exponential weight.

\end{abstract}
\vpn{\small{\bf 2000 Mathematical Subject Classification:} 35Q30, 76D05,
76D07}
\section {Introduction}
\footnotetext[1] {Ri Myong-Hwan: Institute of Mathematics,
Academy of Sciences, DPR Korea, 
email: {\tt mhri@amss.ac.cn}}
\footnotetext[2] {Reinhard Farwig: Department of Mathematics,
Darmstadt University of Technology, 64289 Darmstadt, Germany, 
email:  {\tt farwig@mathematik.tu-darmstadt.de}}

Let
\begin{equation}
\label{E1.2}
\Om=\bigcup_{i=0}^{m}\Om_i
\end{equation}
be a cylindrical domain of $C^{1,1}$-class where $\Om_0$ is a
bounded domain and $\Om_i, i=1,\ldots, m,$ are disjoint semi-infinite
straight cylinders, that is, in possibly different coordinates,
$$ \Om_i=\{x^i=(x^i_1,\ldots, x^i_n)\in\R^n: \:x^i_n>0,\,
   (x^i_1,\ldots, x^i_{n-1})\in\Si^i\},$$
where the cross sections $\Si^i\subset \R^{n-1}, i=1,\ldots,m,$ are bounded domains
and $\Om_i\cap\Om_j=\emptyset$ for $i\neq j$.

Given $\be_i\geq 0, i=1,\ldots,m,$ introduce the space
$$\begin{array}{l}
L^q_\tb(\Om)=\{U\in L^q(\Om): e^{\be_i x^i_n}U|_{\Om_i}\in L^q(\Om_i)\},\ek
\|U\|_{L^q_\tb(\Om)}:=\big(\|U\|^q_{L^q(\Om_0)}
+\sum_{i=1}^m \|e^{\be_i x^i_n}U\|_{L^q(\Om_i)}^q \big)^{1/q}
\end{array}$$
for $1<q<\infty$.
Moreover, let $W^{k,q}_\tb(\Om)$, $k\in\N$, be the space of functions
whose derivatives up to $k$-th order belong to $L^q_\tb(\Om)$,
where a norm is endowed in the standard way.
Let $W^{1,q}_{0,\tb}(\Om)=\{u\in W^{1,q}_\tb(\Om):u|_{\pa\Om}=0\}$.
Let $L^q_\si(\Om)$ and $L^q_{\tb,\si}(\Om)$ be the completion of the set
$C^\infty_{0,\sigma}(\Om)=\{u \in C^\infty_0 (\Om)^n:\, \div u=0\}$
in the norm of $L^q(\Om)$ and $L^q_\tb(\Om)$, respectively.
Then we consider the Stokes operator
$A=A_{q,\tb}=-P_q\Da$ in $L^q_{\tb,\si}(\Om)$ with domain
$${\cal D}(A)=W^{2,q}_\tb(\Om)^n\cap W^{1,q}_{0,\tb}(\Om)^n\cap L^q_\si(\Om),$$
where $P_q$ is the Helmholtz projection of $L^q(\Om)$
onto $L^q_\si(\Om)$.

The goal of this paper is to study resolvent estimates and maximal $L^p$-regularity
of the Stokes operator in Lebesgue spaces with exponential weights in the axial direction.
The semigroup approach to instationary Navier-Stokes equations
is a very convenient tool to prove existence, uniqueness and stability
 of solutions; to this end, resolvent estimates of the Stokes operator
 must be obtained. Moreover, maximal regularity of the Stokes operator helps to deal with the nonlinearity of the Navier-Stokes equations.

There are many papers dealing with resolvent estimates
(\cite{FS94}, \cite{FS97}, \cite{Fr03},
\cite{Fr01}, \cite{Gi81}; see Introduction of \cite{FR05-1}
for more details) or maximal regularity
(see e.g. \cite{Ab04-1}, \cite{Fr02}, \cite{Fr01}) of Stokes operators
for domains with compact as well as noncompact boundaries.
General unbounded domains are considered in \cite{FKS} by replacing
the space $L^q$ by $L^q\cap L^2$ or $L^q+L^2$.
For resolvent estimates and maximal regularity in
unbounded cylinders without exponential weights in the axial direction
we refer the reader e.g. to \cite{FR05-1}-\cite{FR05-4} and \cite{RY}.
For partial results in the Bloch space of
uniformly square integrable functions on a cylinder we refer to \cite{Sc98}.

Further results on stationary Stokes
 and instationary Stokes and Navier-Stokes
systems in unbounded cylindrical domains can be found e.g. in \cite{Am77}, \cite{Be05},
\cite{Ga94-1}, \cite{Ga94-2}, \cite{KaPi12}-\cite{PiZa08},
\cite{Sc98}-\cite{Sp99}.

Despite of some references showing the existence of stationary
 flows in $L^q$-setting (e.g. \cite{NST99-1}, \cite{NST99-2}, \cite{PT98})
and instationary flows in $L^2$-setting (e.g. \cite{Pi05}, \cite{PiZa08})
that converge at $|x|\ra \infty$ to some limit states (Poiseuille flow or zero flow)
 in unbounded cylinders,
resolvent estimates and maximal regularity of the Stokes operator in
$L^q$-spaces with exponential weights on unbounded cylinders
do not seem to have been obtained yet.
\par\bigskip
We start our work with consideration of the Stokes
operator in straight cylinders;
we get resolvent estimate and maximal regularity of
the Stokes operator even in $L^q_\be(\R;L^r_\om(\Si))$, $1<q,r<\infty$, with
exponential weight $e^{\be x_n}$, $\be>0$,
and arbitrary Muckenhoupt weight $\om\in A_r(\R^{n-1})$ with respect to
$x'\in \Si$ (see Section 2 for the definition).
We note that our resolvent estimate gives, in particular when $\la=0$, a new result on the
existence of a unique flow with zero flux for the
 stationary Stokes system in $L^q_\be(\R,L^r_\om(\Si))$.
Next,  based on the results for
straight cylinders, we get resolvent estimates and 
maximal $L^p$-regularity of the Stokes operator in $L^q_\tb(\Om)$,
$1<q<\infty$, for general cylinders $\Om$ using a cut-off technique.

The proofs for straight cylinders are mainly based on the theory of Fourier
analysis. 
By the application of the partial Fourier transform along
the axis of the cylinder $\Si\ti\R$ the {\it generalized Stokes resolvent
system}
$$ \begin{array}{rcl}
\lambda U - \Da U +  \na P & =& F \quad \mbox{ in } \Si\ti\R,\ek
(R_\lambda) \hspace*{4cm} \div U & = & G \quad \mbox{ in }\Si\ti\R,\ek
u & = & 0 \quad \mbox{ on } \pa \Si\ti\R,
\end{array}$$
is reduced to the
{\it parametrized Stokes system} in the cross-section $\Sigma$:
$$ \begin{array}{rcll}
 (\lambda+ {\eta}^2 - \Da') \hat{U}' +  \na' \hat{P} & =& \hat{F}' \quad
 & \mbox{ in } \Sigma,\ek
(\lambda+ {\eta}^2 - \Da') \hat{U}_n +  i\eta \hat{P} & =& \hat{F}_n \quad
& \mbox{ in } \Sigma,\ek
(R_{\lambda,\eta}) \hspace{3cm} \div' \hat{U}'+i\eta\hat{U}_n & = &
   \hat{G} \quad &\mbox{  in }\Sigma,\ek
\hat{U}' = 0, \quad \hat{U}_n &=& 0 \quad &\mbox{ on } \pa   \Sigma,
\end{array}$$
which involves the Fourier phase variable $\eta \in {\mathbb C}$ as
parameter. Now, for fixed $\be\geq 0$ let
$$(\hat{u}, \hat{p}, \hat{f}, \hat{g})(\xi):=(\hat{U}, \hat{P}, \hat{F}, \hat{G})(\xi+i\be).$$
Then $(R_{\la,\eta})$ is reduced to the system
$$ \begin{array}{rcll}
 (\lambda+ (\xi+i\be)^2 - \Da') \hat{u}'(\xi) +  \na' \hat{p}(\xi) & =& \hat{f}'(\xi) \quad
 & \mbox{ in } \Sigma,\ek
(\lambda+ (\xi+i\be)^2 - \Da') \hat{u}_n(\xi) +  i(\xi+i\be) \hat{p}(\xi) & =& \hat{f}_n(\xi) \quad
& \mbox{ in } \Sigma,\ek
(R_{\lambda,\xi,\beta}) \hspace{3cm} \div' \hat{u}'(\xi)+i(\xi+i\be)\hat{u}_n(\xi) & = &
   \hat{g}(\xi) \quad &\mbox{  in }\Sigma,\ek
\hat{u}'(\xi) = 0, \quad \hat{u}_n(\xi) &=& 0 \quad &\mbox{ on } \pa   \Sigma.
\end{array}$$
We will get estimates of solutions
to $(R_{\la,\xi,\be})$ independent of $\xi\in\R^*:=\R\setminus \{0\}$ and
$\la$ in $L^r$-spaces
with Muckenhoupt weights, which yield $\cal{R}$-boundedness of the family
of solution operators $a(\xi)$ for $(R_{\la,\xi,\be})$ with $g=0$ due to
an extrapolation property of operators defined on $L^r$-spaces with
Muckenhoupt weights, see Theorem \ref{T4.8}.
Then, an operator-valued
Fourier multiplier theorem (\cite{We01}) implies the
estimate of $e^{\be x_n}U=\finv (a(\xi) \fou f)$  for the solution $U$ to
$(R_\lambda)$ with $G=0$ in the straight cylinder $\Si\ti\R$.
In order to prove maximal regularity of the Stokes operator in straight cylinders
we use that maximal
regularity of an operator $A$ in a {\it UMD} space $X$ is implied by the
$\cal{R}$-boundedness of the operator family
\begin{equation}
\label{E1.1}
\{\la(\la+A)^{-1}:\,\, \la\in i\,\R\}
\end{equation}
in ${\cal L}(X)$, see \cite{We01}. We show the $\cal{R}$-boundedness of
\eq{1.1} for the Stokes operator $A:=A_{q,r;\be,\om}$ in
$L^q_\be(\R:L^r_\om(\Si))$ by virtue of Schauder decomposition
techniques; to be more precise, we use the Schauder decomposition
$\{\Da_j\}_{j\in\Z}$ where $\Da_j=\finv \chi_{[2^j,2^{j+1})}\fou$
to get $R$-boundedness of the family \eq{1.1}.

The proof for general cylinders, Theorem \ref{T2.4} and Theorem \ref{T2.5}, uses a cut-off technique
based on the result for resolvent estimates and maximal regularity 
without exponential weights in \cite{FR05-4} and
the result (Theorem \ref{T2.3}) for straight cylinders.

\vspace{0.3cm}

This paper is organized as follows. In Section 2 the main results of
this paper (Theorem \ref{T2.1}, Corollary \ref{C2.2}, Theorem \ref{T2.3} --
 Theorem \ref{T2.5}) and
preliminaries are given. In Section 3 we obtain the estimate for $(R_{\lambda,\xi,\be})$ on
bounded domains, see Theorem \ref{T3.8}. In Section 4 proofs of the main
results are given.
%
%
\section{Main Results and Preliminaries}
Let $\Sigma \ti {\mathbb R}$ be an infinite cylinder of ${\mathbb R}^n$
with bounded cross section $\Si\subset \R^{n-1}$ and with generic point
$x \in \Si\ti\R$ written in the form $x = (x', x_n) \in \Si\ti\R$, where $x' \in \Sigma$
and $x_n\in{\mathbb R}$. Similarly, differential operators in ${\mathbb R}^n$
are split, in particular, $\Da=\Da'+\pa^2_n$ and $\na=(\na', \pa_n)$.

For $q\in (1, \infty)$ we use the standard notation
$L^q(\Si\ti\R)=L^q(\R;L^q(\Si))$ for classical Lebesgue spaces with norm
$\|\cdot\|_q=\|\cdot\|_{q;\Si\ti\R}$ and $W^{k,q}(\Si\ti\R), k\in\N,$ for the
usual Sobolev spaces with norm $\|\cdot \|_{k, q;\Si\ti\R}$. We do not
distinguish between spaces of scalar functions and vector-valued
functions as long as no confusion arises. In particular, we use the
short notation $\|u, v\|_X$ for $\|u\|_X+\|v\|_X$, even if $u$ and
$v$ are tensors of different order.

Let $1<r<\infty$. A function $0\leq \om\in L^1_{\text{loc}}(\R^{n-1})$ is called
$A_r$-{\it weight} ({\it Muckenhoupt weight}) on $\R^{n-1}$ iff
$$ {\cal A}_r(\om):=\sup_Q \left(\frac{1}{|Q|}\int_Q\om\, dx'\right)
   \cdot\left(\frac{1}{|Q|}\int_Q\om^{{-1}/{(r-1)}}\, dx'\right)^{r-1}<\infty $$
where the supremum is taken over all cubes of $\R^{n-1}$ and $|Q|$ denotes the
$(n-1)$-dimensional Lebesgue measure of $Q$. We call ${\cal A}_r(\om)$ the
$A_r$-constant of $\om$ and denote the set of all $A_r$-weights on
$\R^{n-1}$ by $A_r=A_r(\R^{n-1})$. Note that
$$ \om\in A_r \quad \text{iff} \quad \om':=
   \om^{-1/{(r-1)}}\in A_{r'},\quad r'=r/{(r-1)},$$
and $A_{r'}(\om')=A_r(\om)^{r'/r}$.
A constant $C=C(\om)$ is called $A_r$-{\it consistent}
if for every $d>0$
$$\sup\,\{C(\om):\; \om\in A_r,\, {\cal A}_r(\om)<d\}<\infty.$$
We write $\om(Q)$ for $\int_Q \om\, dx'$.

Typical Muckenhoupt weights are the radial functions
$\om(x)=|x|^\al$: it is well-known that $\om \in A_r(\R^{n-1})$
if and only if $-(n-1)<\al<(r-1)(n-1)$;
the same bounds for $\al$ hold when $\om(x)=(1+|x|)^\al$ and
$\om(x)=|x|^\al(\log(e+|x|)^\beta$ for all $\beta\in\R.$
For further examples we refer to \cite{FS97}.

Given $\om\in A_r, r\in(1,\infty)$, and an arbitrary domain
$\Si\subset\R^{n-1}$ let
$$ L^r_\om(\Si)=\Big\{u\in L^1_{\text{loc}}(\bar{\Si}): \|u\|_{r,\om} =
   \|u\|_{r,\om;\Si} = \Big(\int_\Si |u|^r\om\,dx'\Big)^{1/r}<\infty\Big\}.$$
For short we will write $L^r_\om$ for $L^r_\om(\Si)$ provided that the
underlying domain $\Si$ is known from the context.
It is well-known that $L^r_\om$ is a separable reflexive Banach space with
dense subspace $C^\infty_0(\Si)$. In particular
$(L^r_\om)^*=L^{r'}_{\om'}$. As usual, $W^{k,r}_\om(\Si)$,
$k\in\N$, denotes the weighted Sobolev space with norm
$$ \|u\|_{k,r,\om}=\Big(\sum_{|\alpha|\leq k}\|D^\alpha u\|^r_{r,\om}\Big)^{1/r},$$
where $|\alpha|=\alpha_1+\cdots+\alpha_{n-1}$ is the length of the
multi-index $\alpha=(\alpha_1,\ldots,\alpha_{n-1})\in\N^{n-1}_0$ and
$D^\alpha=\pa_1^{\alpha_1}\cdot\ldots\cdot \pa_{n-1}^{\alpha_{n-1}}$;
moreover,  $W^{k,r}_{0,\om}(\Si):=\overline{C^\infty_0(\Si)}^{\|\cdot\|_{k,r,\om}}$
and $W^{-k,r}_{0,\om}(\Si):=(W^{k,r'}_{0,\om'}(\Si))^*$, where $r'=r/{(r-1)}$.
We introduce the weighted homogeneous Sobolev space
$$ \hw^{1,r}_\om(\Si)=\left\{u\in {L^1_{\text{loc}}(\bar{\Si})}/\R:\:
   \na' u \in L^r_{\om}(\Si)\right\}$$
with norm $\|\na' u\|_{r,\om}$ and its dual space $\hw^{-1,r'}_{\om'}:=(\hw^{1,r}_\om)^*$
with norm $\|\cdot\|_{-1,r',\om'}= \|\cdot\|_{-1,r',\om';\Si}$.

Let $q,r\in (1,\infty)$. On an infinite cylinder $\Si\ti\R$, where $\Si$ is a
bounded $C^{1,1}$-domain of $\R^{n-1}$, we introduce the function space
$L^q(L^r_\om):=L^q(\R;L^r_\om(\Si))$ with norm
$$ \|u\|_{L^q(L^r_\om)}=\left(\int_\R \Big(\int_\Si |u(x',x_n)|^r
   \om(x')\,dx'\Big)^{q/r}\,dx_n\right)^{1/q}. $$
Furthermore, $W^{k;q,r}_\om(\Si\ti\R), k\in\N,$ denotes the Banach space
of all functions in $\Si\ti\R$ whose derivatives of order up to $k$
belong to $L^q(L^r_\om)$ with norm
$\|u\|_{W^{k;q,r}_\om}=(\sum_{|\alpha| \leq k}
\|D^{\alpha}u\|^2_{L^q(L^r_\om)})^{1/2}$, where $\alpha\in\N^n_0$,
and let $W^{1;q,r}_{0,\om}(\Om)$ be the completion of the set
$C^\infty_0 (\Om)$ in $W^{1;q,r}_\om(\Om)$.
Given $\be>0$, we denote by
$$L^q_\be(L^r_\om):=\{u: e^{\be x_n}u\in L^q(L^r_\om)\}$$
with norm $\|e^{\be x_n}\cdot\|_{L^q(L^r_\om)}$ and for $k\in\N$
$$W^{k;q,r}_{\be,\om}(\Si\ti\R):=\{u: e^{\be x_n}u\in W^{k;q,r}_{\om}(\Si\ti\R)\}$$
with norm $\|e^{\be x_n}\cdot\|_{W^{k;q,r}_{\om}(\Si\ti\R)}$.
Finally, $L^q(L^r_\om)_\sigma$ and $L^q_\be(L^r_\om)_\sigma$ are completions
in the space $L^q(L^r_\om)$ and $L^q_\be(L^r_\om)$ of the set
$$ C^\infty_{0,\sigma}(\Si\ti\R)=\{u \in C^\infty_0 (\Si\ti\R)^n;\quad \div u=0\},$$
respectively.

The Fourier transform in the variable $x_n$ is denoted by $\cal F$
or $\,\widehat{}\,$ and the inverse Fourier transform by ${\cal F}^{-1}$
or $^\vee$.  For $\varepsilon \in (0, \frac{\pi}{2})$ we define the
complex sector
$$ S_\varepsilon = \{ \lambda \in \C; \lambda\neq0, | \text{arg} \lambda |
   < \frac{\pi}{2}+ \varepsilon \}.$$

The first main theorem of this paper is as follows.
%
%
\begin{theo} {\bf (Weighted Resolvent Estimates)}
\label{T2.1}
Let $\Si$ be a bounded domain of $C^{1,1}$-class with $\alpha_0>0$ and $\alpha_1>0$
being the least eigenvalue of the Dirichlet and Neumann Laplacian in $\Sigma$, and let
$\bar\al:=\min\{\al_0,\al_1\}$, $\be\in (0,\sqrt{\bar\al})$,  $\al\in (0, \bar\al-\be^2)$,
$0 <\varepsilon<\ve^*:=\arctan\big({\frac{1}{\be} \sqrt{\bar\al-\be^2-\al}}\big)$, $ 1<q,r<\infty$ and $\om\in A_r.$
Then for every $f \in L^q_\be({\mathbb R}; L^r_\om(\Sigma))$,
and $\lambda\in-\alpha+S_\varepsilon$ there exists a unique solution $(u,\na p)$
to $(R_\lambda)$ (with $g=0$)
such that
$$ (\la+\al)u, \na^2u,\na p\in L^q_\be(L^r_\om)$$
and
\begin{equation}
\label{E2.1}
\|(\lambda+\alpha) u, \na^2u, \na p \|_{L^q_\be(L^r_\om)} \leq C \|f\|_{L^q_\be(L^r_\om)}
\end{equation}
with an $A_r$-consistent constant
$C=C(q,r,\alpha,\be,\ve,\Si,{\cal A}_r(\om))$ independent of $\lambda$.
\end{theo}

In particular we obtain from Theorem 2.1 the following corollary on
resolvent estimates of the Stokes operator in the cylinder $\Om$.
%
%
\begin{coro} {\bf (Stokes Semigroup in Straight Cylinders)}
\label{C2.2}
Let $1<q,r<\infty$, $\om\in A_r(\R^{n-1})$ and define the Stokes
operator $A = A_{q,r;\be,\om}$ on $\Si\ti\R$ by
\begin{equation}
\label{E2.2}
D(A)= W^{2;q,r}_{\be,\om} ( \Si\ti\R) \cap W^{1;q,r}_{0,\be,\om} (\Si\ti\R) \cap L^q_\be(L^r_\om)_\sigma
   \subset L^q_\be(L^r_\om)_\sigma,  \;Au= -P_{q,r;\be,\om}\Da u,
\end{equation}
where $P_{q,r;\be,\om}$ is the Helmholtz projection in
$L^q_\be(L^r_\om)$ (see \cite{Fa03}). Then, for every
$\varepsilon \in (0, \ve^*)$ and $\alpha\in(0,\bar\alpha-\be^2)$, $\be\in (0,\sqrt{\bar\al})$,
$-\alpha + S_\varepsilon$ is contained in the resolvent set of $-A$,
and the estimate
\begin{equation}
\label{E2.3}
\|(\lambda + A)^{-1}\|_{{\cal L}(L^q(L^r_\om)_\sigma)}
\leq {\di {\frac{C}{|\lambda + \alpha|}}}, \quad \forall \lambda
\in -\alpha + S_\varepsilon,
\end{equation}
holds with an $A_r$-consistent constant $C=C(\Si,q,r,\al,\be,
\ve,\car)$.

As a consequence, the Stokes operator generates a bounded analytic semigroup
$\{e^{-tA_{q,r;\be,\om}};t\geq 0 \}$ on $L^q_\be(L^r_\om)_\sigma$ satisfying the estimate
\begin{equation}
\label{E2.4}
 \|e^{-tA_{q,r;\be,\om}}\|_{{\cal L}(L^q_\be(L^r_\om)_\sigma)}\leq C\, e^{-\alpha t}
 \quad \forall \al\in(0,\bar\al-\be^2), \forall t>0,
\end{equation}
with a constant $C=C(q,r,\al,\be,\ve,\Si,\car)$.
\end{coro}

The second important result of this paper is the {\it maximal regularity}
of the Stokes operator in an infinite straight cylinder.
\begin{theo} {\bf {(Maximal Regularity in Straight Cylinders)}}
\label{T2.3}
Let $1<p,q,r<\infty$, $\om\in A_r(\R^{n-1})$ and $\be\in (0,\sqrt{\bar\al})$.
Then the Stokes operator $A=A_{q,r;\be,\om}$ has maximal regularity in
$L^q_\be(L^r_\om)_\si$. To be more precise,
for each $F\in L^p(\R_+; L^q_\be(L^r_\om)_\si)$
the instationary problem
\begin{equation}
\label{E2.4n}
U_t+AU=F,\quad U(0)=0,
\end{equation}
has a unique solution
$U\in W^{1,p}(\R_+; L^q_\be(L^r_\om)_\si)\cap L^p(\R_+; D(A))$
such that
\begin{equation}
\label{E2.5b} \|U, U_t, AU\|_{L^p(\R_+; L^q_\be(L^r_\om)_\si)}
   \leq C\|F\|_{L^p(\R_+; L^q_\be(L^r_\om)_\si)}.
\end{equation}
Analogously, for every $F\in L^p(\R_+; L^q_\be(L^r_\om))$,
the instationary system
\begin{equation*}
U_t-\Delta U+\na P=F,\quad\div U=0,\quad U(0)=0,
\end{equation*}
has a unique solution
$$(U,\na P)\in
\big(W^{1,p}(\R_+; L^q_\be(L^r_\om)_\si)\cap L^p(\R_+; D(A))\big)
\times L^p(\R_+; L^q_\be(L^r_\om))$$
 satisfying the a priori estimate
\begin{equation}
\label{E2.5c}
\|U_t, U,\na U, \na^2U, \na P\|_{L^p(\R_+; L^q_\be(L^r_\om))}
   \leq C\|F\|_{L^p(\R_+; L^q_\be(L^r_\om))}
\end{equation}
with $C=C(\Si,q,r, \be, {\cal A}_r(\om))$.
Moreover, if
$e^{\al t} F\in L^p(\R_+; L^q_\be(L^r_\om))$ for some $\al\in(0,\bar\al-\be^2)$,
then the solution $u$ satisfies the estimate
\begin{equation}
\label{E2.5d}  \|e^{\al t}U, e^{\al t}U_t, e^{\al t}\na^2 U\|_{L^p(\R_+; L^q_\be(L^r_\om))}
   \leq C\|e^{\al t}F\|_{L^p(\R_+; L^q_\be(L^r_\om))}
\end{equation}
with $C=C(\Si,q,r,\al,\be,\car)$.
\end{theo}

As a corollary of Theorem \ref{T2.3} we get the maximal regularity
result for general cylinder $\Om$  with several exits to infinity given by \eq{1.2}.

\begin{theo} {\bf (Stokes Semigroup in General Cylinders)}
\label{T2.4}
Let a $C^{1,1}$-domain $\Om$ be given by \eq{1.2} and  $\be_i>0$ for $i=1,\ldots,m$
satisfy the same assumptions on $\be$ with $\Si^i$ in place of $\Si$.
Then, the Stokes operator $A_{q,\tb}(\Om)$ generates an
exponentially decaying analytic semigroup $\{e^{-tA_{q,\tb}}\}_{t\geq 0}$ in $L^q_{\tb,\si}(\Om)$.
 \end{theo}

\begin{theo} {\bf (Maximal Regularity in General Cylinders)}
\label{T2.5}
Let a $C^{1,1}$-domain $\Om$ be given by \eq{1.2} and  $\be_i>0$ for $i=1,\ldots,m$
satisfy the same assumptions on $\be$ with $\Si^i$ in place of $\Si$.
Then, the Stokes operator $A_{q,\tb}$ has maximal regularity in $L^q_{\tb,\si}(\Om)$;
to be more precise, for any  $F\in L^p(\R_+; L^q_{\tb,\si}(\Om))$
the Cauchy problem
\begin{equation}
\label{E2.8}
U_t+A_{q,\tb}U=F,\; U(0)=0, \quad \text{in } L^q_{\tb,\si}(\Om),
\end{equation}
has a unique solution $U$ such that
\begin{equation}
\label{E2.9}
\|U, U_t, A_{q,\tb}U\|_{L^p(\R+;L^q_{\tb,\si}(\Om))}\leq C\|F\|_{L^p(\R+;L^q_{\tb,\si}(\Om))}
\end{equation}
with some constant $C=C(q,\Om)$.

Equivalently, if $F\in L^p(\R_+; L^q_\tb(\Om))$,
then the instationary Stokes system
\begin{equation}
\label{E2.7}
\begin{array}{rcl}
U_t - \Da U +  \na P & =& F \quad \mbox{ in } \R_+\ti\Om,\ek
\div U & = & 0 \quad \mbox{ in }\R_+\ti\Om,\ek
U(0) & = & 0 \quad \mbox{ in }\Om,\ek
U & = & 0 \quad \mbox{ on } \pa \Om,
\end{array}
\end{equation}
has a unique solution $(U,\na P)$ such that
\begin{equation}
\label{E2.6}
\begin{array}{l}
(U,\na P)\in
(L^p(\R_+; W^{2,q}_{\tb}(\Om)\cap W^{1,q}_{0}(\Om))\cap L^q_\si(\Om))
\ti L^p(\R_+; L^q_\tb(\Om)),\ek
 U_t\in L^p(\R_+; L^q_\tb(\Om)),\ek
\|U\|_{L^p(\R_+; W^{2,q}_\tb(\Om)\cap W^{1,q}_{0}(\Om))}+
\|U_t, \na P\|_{L^p(\R_+; L^q_\tb(\Om))}\leq
C\|F\|_{L^p(\R_+; L^q_\tb(\Om))}.
\end{array}
\end{equation}
\end{theo}

\begin{rem}
\label{R2.5} {\rm  We note that in \eq{2.4n} and in \eq{2.7} we may take nonzero
initial values $u(0)=u_0$ in the interpolation space
$(L^q_\be(L^r_\om)_\si, D(A_{q,r;\be,\om}))_{1-1/p,p}$
and
$U(0)=U_0\in (L^q_\tb(\Om),W^{2,q}_\tb(\Om)\cap W^{1,q}_{0,\tb}(\Om))_{1-1/p,p}$,
respectively.

}
\end{rem}

For the proofs in Section 3 and Section 4, we need some preliminary results
for Muckenhoupt weights.
%
%
\begin{propo}{\rm (\cite{Fa03}, Lemma 2.4)}
\label{P2.5}
Let $1<r<\infty$ and $\om\in A_r(\R^{n-1})$.

(1) Let $T: \R^{n-1}\rightarrow  \R^{n-1}$ be a bijective, bi-Lipschitz vector field.
Then, it holds that $\om \circ T\in A_r(\R^{n-1})$ and ${\cal A}_r(\om \circ T)
\leq c\,{\cal A}_r(\om)$ with a constant $c=c(T,r)>0$ independent of $\om$.

(2)  Define the weight $\tilde{\om}(x')=\om(|x_1|,x'')$ for $x'=(x_1,x'')\in \R^{n-1}$.
Then $\tilde{\om}\in A_r$ and ${\cal A}_r(\tilde{\om})\leq 2^r\,{\cal A}_r(\om)$.

(3)  Let $\Si\subset\R^{n-1}$ be a bounded domain. Then there exist
$\tilde{s},s\in (1,\infty)$ satisfying
$$L^{\tilde{s}}(\Si)\hookrightarrow L^r_\om(\Si)\hookrightarrow L^s(\Si).$$
Here $\tilde{s}$ and $\frac{1}{s}$ are $A_r$-consistent. Moreover,
the embedding constants can be chosen uniformly on a set $W\subset A_r$ provided that
\begin{equation}
\label{E2.5}
\sup_{\om\in W}{\cal A}_r(\om)<\infty,\quad \int_Q\om\, dx'=1\quad
\text{for all  }\,\om\in W,
\end{equation}
for a cube $Q\subset \R^{n-1}$ with $\bar{\Si}\subset Q$.
\end{propo}
%
%
\begin{propo}{\rm (\cite{Fa03}, Proposition 2.5)}
\label{P2.6} Let $\Si\subset\R^{n-1}$ be a bounded Lipschitz domain
and let $1<r<\infty$.

(1) For every $\om\in A_r$ the continuous embedding
$W^{1,r}_\om(\Si)\hookrightarrow L^r_\om(\Si)$ is compact.

(2) Consider a sequence of weights $(\om_j)\subset A_r$ satisfying \eq{2.5} for
$W=\{\om_j:j\in\N\}$ and a fixed cube $Q\subset\R^{n-1}$ with $\bar{\Si}\subset Q$.
Further let $(u_j)$ be a sequence of functions on $\Si$ satisfying
$$  \sup_j \|u_j\|_{1,r,\om_j}<\infty \quad \text{and}\quad u_j\wrar 0\quad
    \text{in }W^{1,s}(\Si)$$
for $j\rightarrow \infty$ where $s$ is given by Proposition \ref{P2.5} (3). Then
$$ \|u_j\|_{r,\om_j}\rightarrow 0 \quad \text{for   }j\rightarrow \vspace{0.2cm} \infty.$$

(3) Under the same assumptions on $(\om_j)\subset A_r$ as in (2) consider a
sequence of functions $(v_j)$ on $\Si$ satisfying
$$ \sup_j \|v_j\|_{r,\om_j}<\infty \quad \text{and}\quad v_j\wrar 0
   \quad \text{in }L^s(\Si)$$
for $j\rightarrow \infty$. Then considering $v_j$ as functionals on
$W^{1,r'}_{\om'_j}(\Si)$
$$ \|v_j\|_{(W^{1,r'}_{\om'_j}(\Si))^*}\rightarrow 0
   \quad \text{for   }j\rightarrow \infty.$$
\end{propo}
%
%
\begin{propo}
\label{P2.7}
Let $r\in(1,\infty),\; \om\in A_r$ and
$\Si\subset\R^{n-1}$ be a bounded Lipschitz domain. Then there exists
an $A_r$-consistent constant $c=c(r,\Si,{\cal A}_r(\om))>0$ such that
$$\|u\|_{r,\om}\leq c\|\na'u\|_{r,\om}$$
for all $u\in W^{1,r}_\om(\Si)$ with vanishing integral mean $\int_\Si u\,dx'=0$.
\end{propo}
{\bf Proof:} See the proof of \cite{Fr01}, Corollary 2.1 and its conclusions;
checking the proof, one sees that the constant $c=c(r,\Si, {\cal A}_r(\om))$ is
$A_r$-consistent.
\hfill\qed
\medskip

Finally we cite the Fourier multiplier theorem in weighted spaces.
%
%
\begin{theo}  \label{T2.6} {\rm (\cite{GR85}, Ch. IV, Theorem 3.9)}
Let $m\in C^k(\R^k\setminus\{0\}),k\in\N,$ admit a constant $M\in \R$ such that
$$ |\eta|^\gamma |D^\gamma m(\eta)|\leq M
   \quad \text{for all}\quad \eta\in \R^k\setminus\{0\}$$
and multi-indices $\gamma\in\N^k_0$ with $|\gamma|\leq k$.
Then for all $1<r<\infty$ and $\om\in A_r(\R^k)$ the multiplier operator
$Tf={\cal F}^{-1}m(\cdot){\cal F}f$ defined for all rapidly decreasing functions
$f\in{\cal S}(\R^k)$ can be uniquely extended to a bounded linear operator from
$L^r_\om(\R^k)$ to $L^r_\om(\R^k)$. Moreover, there exists an $A_r$-consistent
constant $C=C(r,\car)$ such that
$$\|Tf\|_{r,\om}\leq CM\|f\|_{r,\om}\,,\quad f\in L^r_\om(\R^k)\,.$$
\end{theo}
%
%
%
%
%
\section{Resolvent estimate of the Stokes operator in weighted spaces on infinite straight cylinders}
In this section we  obtain the resolvent estimate of the Stokes operator
in Lebesgue spaces with exponential weight with respect to the axial variable and Muckenhoupt weight
for cross-sectional variables in an infinite straight cylinder $\Si\ti\R$, where the cross-section $\Si$  is a $C^{1,1}$-bounded domain.

\subsection {Estimate for the problem $(R_{\la,\xi,\be})$}

In this subsection we get estimates for $(R_{\la,\xi,\be})$ 
independent of $\la$ and $\xi\in\R^*$ in $L^r$-spaces with Muckenhoupt weights. To this aim we rely partly on cut-off techniques using the
results for $(R_{\la,\xi})$ (i.e., the case $\be=0$)
in the whole and bent half spaces in \cite{FR05-3} (Theorem \ref{T3.1} below).
 The main existence and uniqueness result
in weighted $L^r$-spaces for $(R_{\la,\xi,\be})$
is described in Theorem \ref{T3.8}.

For whole or bent half spaces $\Si$, $g\in \hw^{-1,r}_\om(\Si)+L^r_{\om}(\Si)$
and $\eta=\xi+i\be, \xi\in \R^*, \be\geq 0$,
we use notation
$$ \|g; \hw^{-1,r}_\om+L^r_{\om,1/\eta}\|=\inf
\{\|g_0\|_{-1,r,\om}+\|g_1/\eta\|_{r,\om}: g=g_0+g_1, g_0\in \hw^{-1,r}_\om, g_1\in L^r_{\om}\}.$$
In the following we put $R_{\la,\xi}\equiv R_{\la,\xi,0}$.

\begin{theo}
\label{T3.1}
 Let $n \geq 3,\, 1<r<\infty,\, \om\in A_r(\R^{n-1})$,
$0<\ve <\frac {\pi}{2}$, $\xi \in {\mathbb R}^*$,$\la \in S_\ve$, $0< \ve <\pi /2$
and $\mu =|\la +\xi^2|^{1/2}$.
\begin{itemize}
\item[(i)] (\cite{FR05-3}, Theorem 3.1)
Let $\Si={\mathbb R}^{n-1}$.
If $f \in L^r_\om(\Si)$ and $ g \in W^{1,r}_\om(\Si)$, then the problem
$(R_{\la, \xi})$ has a unique solution $(u,p) \in W^{2,r}_\om(\Si) \ti W^{1,r}_\om(\Si)$
satisfying
\begin{equation}
\label{E3.9}
  \| \mu ^2u, \mu \na 'u, \na '^2u, \na 'p, \xi p\|_{r,\om}
 \leq c\big(\|f, \na 'g, \xi g\|_{r,\om}+\| \la g; \hw^{-1,r}_\om+L^r_{\om,1/\xi}\|\big)
\end{equation}
with an $A_r$-consistent constant $c=c(\ve,r,{\cal A}_r(\om))$  independent of
$\la$ and $\xi$.

\item[(ii)] (\cite{FR05-3}, Theorem 3.5) Let
$$\Si=H_\si=\{x'=(x_1,x'');\, x_1>\si(x''),x''\in {\mathbb R}^{n-2}\}$$
for a given function $\si \in C^{1,1}({\mathbb R}^{n-2})$.
Then there are $A_r$-consistent constants
$K_0=K_0(r,\ve,{\cal A}_r(\om))>0$ and $\la_0=\la_0(r,\ve,{\cal A}_r(\om))>0$ independent of
$\la$ and $\xi$ such that,
if $\| \na'\si \|_\infty \leq K_0$, for every
$f\in L^r_\om(\Si)$ and $ g\in W^{1,r}_\om(\Si)$
the problem
$(R_{\la,\xi})$ has a unique solution
$(u,p)\in (W^{2,r}_\om(\Si)\cap W^{1,r}_{0,\om}(\Si))\ti W^{1,r}_\om(\Si)$.
This solution satisfies the estimate
\begin{equation}
\label{E3.31}
\begin{array}{l}
\|\mu^2u,\mu\na'u,\na'^2u,\na'p,\xi p\|_{r,\om}\ek
\hspace{3cm}\leq c\big(\|f,\na'g,\xi g\|_{r,\om}+
\|\la g;\hw^{-1,r}_\om(\Si)+L^r_{\om,1/\xi}(\Si)\|\big)
\end{array}
\end{equation}
with an $A_r$-consistent constant $c=c(r,\ve,{\cal A}_r(\om))$.
\end{itemize}

\end{theo}

On the bounded domain $\Si\subset \R^{n-1}$ of $C^{1,1}$-class let $\alpha_0$
and $\al_1$
denote the smallest eigenvalue of the Dirichlet and Neumann Laplacian, respectively,
 i.e.,
\begin{equation}
\label{E3.10}
\begin{array}{l}
\alpha_0:=\inf \{\|\na u\|^2_2:\: u\in W^{1,2}_0(\Si), \|u\|_2=1\}>0,
\ek
\alpha_1:=\inf \{\|\na u\|^2_2:\: u\in W^{1,2}(\Si), \frac{\pa u}{\pa n}|_{\pa\Si}=0,\|u\|_2=1\}>0,\ek
\bar\al:=\min\{\al_0,\al_1\}.
\end{array}
\end{equation}

For fixed $\la\in \C\setminus (-\infty,-\alpha_0]$, $\eta=\xi+i\be$,
$\xi\in \R^*$, $\be\geq 0$,
 and $\om\in A_r$
we introduce the {\it parametrized Stokes operator} $S=S^\om_{r,\la,\eta}$ by
$$ S(u,p)= \left( \begin{array}{c}
           (\la+\eta^2-\Da')u'+\na'p\\
           (\la+\eta^2-\Da')u_n+i\eta p\\
           -\,\div\hspace{-0.1cm}_\eta u
\end{array} \right)$$
defined on ${\cal D}(S)={\cal D}(\Da'_{r,\om})\ti W^{1,r}_\om(\Si)$, where
${\cal D}(\Da'_{r,\om})=W^{2,r}_\om(\Si)\cap W^{1,r}_{0,\om}(\Si)$ and
$$\div\hspace{-0.1cm}_\eta u =\div'u'+i\eta u_n.$$
For $\om\equiv 1$ the operator $S^\om_{r,\la,\eta}$ will be denoted by $S_{r,\la,\eta}$.
Note that the image of ${\cal D}(S)$ by $\div\hspace{-0.1cm}_{\eta}$
is included in $W^{1,r}_\om(\Si)$ and
$W^{1,r}_\om(\Si)\subset L^r_{0,\om}(\Si)+L^r_{\om}(\Si)$, where
$$ L^r_{0,\om}(\Si):=\big\{u\in L^r_\om(\Si):\: \int_{\Si}u\,dx'=0 \big\}.$$
Using Poincar\'{e}'s inequality in weighted spaces, see Proposition \ref{P2.7},
one can easily check the continuous embedding
$L^r_{0,\om}(\Si)\hookrightarrow \hw^{-1,r}_\om(\Si)$; more precisely,
$$\|u\|_{-1,r,\om}\leq c\|u\|_{r,\om}\,, \quad u\in L^r_{0,\om}(\Si),$$
with an $A_r$-consistent constant $c>0$.   For bounded domain $\Si$ we use the notation
$$ \|g;L^r_{0,\om}+L^r_{\om,1/\eta}\|_0:=
   \inf \{\|g_0\|_{-1,r,\om}+\|g_1/\eta\|_{r,\om}:\:
   g=g_0+g_1, g_0 \in L^r_{0,\om}, g_1\in L^r_\om \};$$
note that this norm is equivalent to the norm $\|\cdot\|_{(W^{1,r'}_{\om',\eta})^*}$
where $W^{1,r'}_{\om',\eta}$ is the usual weighted Sobolev space on $\Si$ with norm
$\|\na' u,\eta u\|_{r',\om'}$.

First we consider Hilbert space setting of $(R_{\la,\xi,\be})$.
For $\eta=\xi+i\be$, $\xi\in \R^*, \beta\geq 0$, let us introduce
a closed subspace of $W^{1,r}_0(\Si)$ as
$$V_\eta:=\{u\in W^{1,r}_0(\Si):\div_\eta u=0\}.$$

\begin{lem}
\label{L3.1}
Let $\phi=(\phi',\phi_n)\in W^{-1,2}(\Si)$ be such that
$\lan \phi, v\ran_{W^{-1,2}(\Si), W^{1,2}_0(\Si)}=0$ for all $v\in W^{1,2}_0(\Si)$.
Then, there is some $p\in L^2(\Si)$ with $\phi=(\na p, i\eta p)$.
\end{lem}
{\bf Proof:} This lemma can be proved just by copy of the proof of \cite{FR05-1}, Lemma 3.1
with $\xi\in\R^*$ replaced by $\eta=\xi+i\be$.\qed

\begin{lem}
\label{L3.2}
 \begin{itemize}
 \item[(i)] For any $g\in W^{1,2}(\Si)$, $\eta=\xi+i\be$, $x\in\R^*, \be\geq 0$,
 the equation $\div_\eta u=g$ has at least one solution $u\in W^{2,2}(\Si)\cap W^{1,2}_0(\Si)$
 and
 $$\|u\|_{2,2}\leq c(\|g\|_{1,2}+\frac{1}{\eta}\int_\Si g\,dx'),$$
where $c$ is independent of $g$.
\item[(ii)] Let $\ve\in (0,\pi/2)$, $\be\in (0,\sqrt{\al_0})$ and
\begin{equation}
\label{E3.2}
\la\in \{-\al_0+\be^2+S_\ve\} \cap \{\la\in\C:\Re \la>-\frac{(\Im\la)^2}{4\be^2}-\al_0+\be^2\}.
\end{equation}
Then,
for any $f\in L^2(\Si)$, $g\in W^{1,2}(\Si)$
the system $(R_{\la,\xi,\be})$ has a unique solution $(u,p)\in (W^{2,2}(\Si)\cap W^{1,2}(\Si))\ti W^{1,2}(\Si)$.
\end{itemize}
\end{lem}
{\bf Proof:} -- {\it Proof of (i)}: Let a scalar function $w\in C^\infty_0(\Si)$ be such that
$\int_\Si w\,dx'=0$. Given $g\in W^{1,2}(\Si)$, let $\bar{g}=\int_\Si g\,dx$ and consider a
divergence problem in $\Si$, that is,
$$\div'u'=g-\bar{g}w, \quad u'|_{\pa\Si}=0,$$
which has a solution $u'\in W^{2,2}(\Si)\cap W^{1,2}_0(\Si)$ with
$\|u'\|_{2,2}\leq c\|\na(g-\bar{g}w)\|_2\leq c\|g\|_{1,2}$, by \cite{FS94}, Theorem 1.2.
Then, $u:=(u',\frac{\bar{g}w}{i\eta})$ satisfies $\div_\eta u=g$ and required estimate.

-- {\it Proof of (ii)}: By the assertion (i) of the lemma, we may assume w.l.o.g. that
 $g\equiv 0$. Now, for fixed $\la\in -\al_0+\be^2+S_\ve$ define the
 bilinear form $b: V_\eta\ti V_\eta\mapsto \R$ by
 $$b(u,v):=\int_\Si \big((\la+\eta^2)u\cdot \bar{v}+\na'u\cdot\na'\bar{v}\big)\,dx'.$$
Obviously, $b$ is continuous in $V_\eta\ti V_\eta$. Moreover, $b$ is coercive, that is,
\begin{equation}
\label{E3.1}
|b(u,u)|\geq l(\la,\xi,\be)\|u\|^2_{1,2}
\end{equation}
with some $l(\la,\xi,\be)>0$.
In fact,
$$b(u,u)=\int_\Si ((\Re\la+\xi^2-\be^2)|u|^2+|\na' u|^2)\,dx'
  +i\int_\Si (\Im\la+2\xi\be)|u|^2\,dx'.$$
Note that, due to Poincar\'e's inequality,
$$(\xi^2-\al_0)\|u\|_2^2+\|\na'u\|_2^2>0, \,\forall \xi\in\R^*.$$
Hence, if $\Re\la+\al_0-\be^2\geq 0$, then
$$|b(u,u)|\geq |\int_\Si ((\Re\la+\xi^2-\be^2)|u|^2+|\na' u|^2)\,dx'|
\geq (\xi^2-\al_0)\|u\|_2^2+\|\na'u\|_2^2,$$
where
$$(\xi^2-\al_0)\|u\|_2^2+\|\na'u\|_2^2\geq \|\na'u\|_2^2$$
if $\xi^2-\al_0\geq 0$ and
 $$(\xi^2-\al_0)\|u\|_2^2+\|\na'u\|_2^2\geq (\xi^2/\al_0-1)\|\na'u\|_2^2+\|\na'u\|_2^2\geq \frac{\xi^2}{\al_0}\|\na'u\|_2^2$$
if $\xi^2-\al_0<0$.

Therefore, it remained to prove \eq{3.1} for the case  $\Re\la+\al_0-\be^2< 0$.

Note that if $\Im\la+2\xi\be\neq 0$ then $(R_{\la, \xi,\be})$ coincides
with $(R_{\la_1,\xi})$ where $\la_1=\la-\be^2+2i\xi\be \in -\al_0+S_{\ve_1}$ with
$\ve_1=\max\{\ve, \arctan{\frac{|\Re\la+\al_0-\be^2|}{|\Im\la+2\xi\be|}}\}\in (0,\pi/2)$.
Hence, \eq{3.1} can be proved in the same way as the proof of \cite{FR05-1}, Lemma 3.2 (ii).

Now, suppose that
$$\Im\la+2\xi\be= 0, \quad i.e.,\quad \xi=-\frac{\Im\la}{2\be}.$$
Since \eq{3.1} is trivial for the case $\Re\la+\xi^2-\be^2\geq 0$, we assume
that $$\Re\la+\xi^2-\be^2< 0.$$
In this case, note that due to the condition $\Re\la+\frac{(\Im\la)^2}{4\be^2}-\be^2>-\al_0$
there is some $c(\la,\be)>0$ such that
$$0>\Re\la+\frac{(\Im\la)^2}{4\be^2}-\be^2> c(\la,\be)-\al_0,  \quad c(\la,\be)-\al_0<0.$$
Then,
$$\begin{array}{rcl}
|b(u,u)|&\geq& \int_\Si ((\Re\la+\frac{(\Im\la)^2}{4\be^2}-\be^2)|u|^2+|\na' u|^2)\,dx'\ek
       &\geq&  \int_\Si (c(\la,\be)-\al_0)|u|^2+|\na' u|^2)\,dx'\ek
       &\geq&  \frac{c(\la,\be)}{\al_0}\|\na' u\|_2^2.
\end{array}$$
Finally, \eq{3.1} is proved.

By Lax-Milgram's lemma in view of \eq{3.1}, the variational problem
$$b(u,v)=\int_\Si f \cdot \bar{v}\,dx', \forall v\in V_\eta,$$
has a unique solution $u$ in $V_\eta$. Then, by Lemma \ref{L3.1}, there is some $p\in L^2(\Si)$
such that
$$(\la+\eta^2-\Da')u'+\na' p=f', (\la+\eta^2-\Da')u_n+i\eta p=f_n.$$
Now, applying the well-known regularity theory for Stokes system and Poisson's equation in $\Si$
to
$$-\Da'u'+\na' p=f'-(\la+\eta^2)u', \div'u'=-i\eta u_n, u'|_{\pa\Si}=0$$
and
$$-\Da'u_n=f_n-(\la+\eta^2)u_n-i\eta p, u_n|_{\pa\Si}=0,$$
respectively, we have $(u,p)\in (W^{2,2}(\Si)\cap W^{1,2}_0(\Si))\ti W^{1,2}(\Si)$.
Thus, the assertion (ii) of the lemma is proved. \qed

\begin{rem}
\label{R3.3}
It is seen by elementary calculation that 
the assumption \eq{3.2} on $\la$ of Lemma \ref{L3.2} is satisfied 
for all $\la\in -\al+S_\ve$ if either $\al\in (0, \al_0-\be^2)$ 
and $\ve\in (0, \arctan{\frac{\sqrt{\al_0-\be^2-\al}}{\be}})$ 
or if $\al\in (0, \bar{\al}-\be^2)$ and
 $\ve\in (0, \arctan{\frac{\sqrt{\bar{\al}-\be^2-\al}}{\be}})$. Note that $\bar\al<\al_0$, see \eq{3.10}.

\end{rem}

Now, we turn in considering $(R_{\la,\xi,\be})$ in weighted spaces with weights w.r.t. cross-section as well.
%
%
\begin{lem}\
\label{L3.4}
Let $\xi\in\R^*$, $\be \in (0, \sqrt{\al_0})$, $\al\in (0,\al_0-\be^2)$,
 $\ve\in (0,\arctan{\frac{\sqrt{\al_0-\be^2-\al}}{\be}})$,
$\la\in -\alpha+S_\ve$,
 and $\om\in A_r$, $1<r<\infty$. Then the operator $S=S^\om_{r,\la,\eta}$ is injective and the range
${\cal R}(S)$ of $S$ is dense in $L^r_\om(\Si)\ti W^{1,r}_\om(\Si)$.
\end{lem}
{\bf Proof:} Since, by Proposition \ref{P2.5} (3), there is an $s\in (1,r)$ such that
$L^r_\om(\Si)\subset L^s(\Si)$, one sees immediately that
${\cal D}(S^\om_{r,\la,\eta})\subset {\cal D}(S_{s,\la,\eta}).$
Therefore, $S^\om_{r,\la,\eta}(u,p)=0$ for some $(u,p)\in {\cal D}(S^\om_{r,\la,\eta})$
yields $(u,p)\in {\cal D}(S_{s,\la,\eta})$ and $S_{s,\la,\eta}(u,p)=0$.
Here note that $S_{s,\la,\eta}(u,p)=0$ implies
$$S_{s,\la,\eta}(u,p)=((\be^2-2i\xi\be)u',(\be^2-2i\xi\be)u_n+\be p,\be u_n )^T.$$
Hence, by applying  \cite{FR05-1}, Theorem 3.4 finite number of times and the Sobolev embedding
theorem, we get that $(u,p)\in (W^{2,2}(\Si)\cap W^{1,2}_0(\Si))
\ti W^{1,2}(\Si)$. Therefore, by Lemma \ref{L3.2} we get that $(u,p)=0$, i.e.,
$S^\om_{r,\la,\eta}$ is injective.
\par
On the other hand, by Proposition \ref{P2.5} (3), there is an
$\tilde{s}\in (r,\infty)$ such that
$S_{\tilde{s},\la,\eta}\subset S^\om_{r,\la,\eta}$.
Moreover, by Lemma \ref{L3.2}, for every $(f, g)\in C^\infty_0(\Si)\ti C^\infty(\bar\Si)$,
there is some $(u,p)\in D(S_{2,\la,\eta})$ with $S_{2,\la,\eta}(u,p)=(f,-g)$.
Applying the regularity result of \cite{FS94}, Theorem 1,2 for the Stokes resolvent system in $\Si$
finite number of times  using the Sobolev embedding
theorem,
it can be seen that $(u,p)\in D(S_{q,\la,\eta})$ for all $q\in (1,\infty)$, in particular,
for $q=\tilde{s}$.
Therefore,
$$ C^\infty_0(\Si)\ti C^{\infty}(\bar\Si)
   \subset {\cal R}(S_{\tilde{s},\la,\eta})\subset{\cal R}(S^\om_{r,\la,\eta})
   \subset  L^r_\om(\Si)\ti W^{1,r}_\om(\Si), $$
which proves the assertion on the denseness of ${\cal R}(S)$.
\hfill\qed\vspace{0.3cm}
\par The following lemma gives a preliminary {\it a priori} estimate
for a solution $(u,p)$ of $S(u,p)=(f,-g)$.
%
%
\begin{lem}
\label{L3.5}
Assume the same for $r,\om,\al,\be$ and $\la$ as in Lemma \ref{L3.4}.
Then there exists an $A_r$-consistent constant $c=c(\ve, r, \be, \Si,\car)>0$
such that for
every $(u,p) \in {\cal D}(S_{r,\la,\eta}^\om)$,
\begin{equation}
\label{E4.1}
\begin{array}{l}
\hspace{-0.3cm} \|\mu^2_+u,\mu_+\na'u,\na'^2u,\na'p,\eta p\|_{r,\om}
    \leq c\big(\|f,\na'g,g,\xi g\|_{r,\om}+|\la|\|g;
    L^r_{0,\om}+ L^r_{\om,1/\eta}\|_0\\[1ex]
\hspace{56mm} +\|\na'u,\xi u,p\|_{r,\om} +
|\la|\| u\|_{(W^{1,r'}_{\om'})^*}\big),
\end{array}
\end{equation}
where $\mu_+=|\la+\alpha+\xi^2|^{1/2}, (f,-g)=S(u,p)$ and
$(W^{1,r'}_{\om'})^*$ denotes the dual space of $W^{1,r'}_{\om'}(\Si)$.
\end{lem}
{\bf Proof:} The proof is devided into two parts, i.e, the case $\xi^2>\be^2$
 and the other case $\xi^2\leq \be^2$.

 The proof of the case $\xi^2>\be^2$
is based on a partition of unity in $\Si$ and on the
localization procedure reducing the problem to a finite number of
problems of type
$(R_{\la,\xi})$ in bent half spaces and in the whole space ${\mathbb R}^{n-1}$.
Since $\pa\Si\in C^{1,1}$, we can cover $\pa\Si$ by a
finite number of balls $B_j, j\geq 1$, such that, after a translation
and rotation of coordinates, $\Si\cap B_j$
locally coincides with a bent half space $\Si_j=\Si_{\si_j}$
where $\si_j\in C^{1,1}(\R^{n-1})$ has a compact support,
$\si_j(0)=0$ and $\na'' \si_j(0)=0$. Choosing the balls $B_j$
small enough (and its number large enough) we may assume that
$\|\na''\si_j\|_{\infty}\leq K_0(\ve,r,\Si,\car)$ for all $j\geq 1$ where $K_0$
was introduced in Theorem  \ref{T3.1} (ii).

 According to the covering
$\pa\Si\subset \bigcup_{j\geq 1} B_j$ there are cut-off functions
$(\vp_j)_{j=0}^m$ such that
 such that
\begin{equation}
\label{E3.8}
0\leq \varphi_0, \varphi_j\in C^{\infty}({\mathbb R}^{n-1}),\;
 \sum_{j\geq 0}\varphi_j\equiv 1 \text{ in }\Si,
   \; \text{supp } \varphi_0\subset \Si,\text{supp}\, \varphi_j \subset B_j, j\geq 1.
\end{equation}

Given $(u,p)\in {\cal D}(S)$ and $(f,-g)=S(u,p)$, we get for each
$\varphi_j,\, j\geq 0,$ the local $(R_{\la,\xi})$-problems
\begin{equation}
\label{E4.2}
\begin{array}{rcl}
(\la+\xi^2-\Da')(\varphi_ju')+\na'(\varphi_jp)  & = & f'_j\ek
(\la+\xi^2-\Da')(\varphi_ju_n)+i\xi(\varphi_jp) & = & f_{jn}\ek
       \div\hspace{-0.1cm}_{\xi}(\varphi_ju)    & = & g_j
\end{array}
\end{equation}
for $(\varphi_ju, \varphi_j p), j\geq 0$, in ${\mathbb R}^{n-1}$ or $\Si_j$; here
\begin{equation}
\label{E4.3}
\begin{array}{rcl}
f'_j&=&\varphi_jf'-2\na'\varphi_j\cdot\na'u'-(\Da'\varphi_j)u'
+(\be^2-2i\xi) (\vp_j u')+(\na'\varphi_j)p\ek
f_{jn}&=&\varphi_jf_n-2\na'\varphi_j\cdot\na'u_n-(\Da'\varphi_j)u_n+(\be^2-2i\xi) (\vp_j u_n)+\be(\vp_j p)\ek
g_j&=& \varphi_jg+\na'\varphi_j\cdot u'+\be\vp_j u_n.
\end{array}
\end{equation}

To control $f_j$ and $g_j$ note that $u=0$ on $\pa\Si$;
hence Poincar\'{e}'s inequality
for Muckenhoupt weighted space yields for all $j\geq 0$ the estimate
\begin{equation}
\label{E4.4}
\|f_j,\na'g_j,\xi g_j\|_{r,\om;\Si_j}\leq c(\|f,\na'g,g,\xi g\|_{r,\om;\Si} +
   \|\na'u,\xi u,p\|_{r,\om;\Si}),
\end{equation}
where $\Si_0\equiv \R^{n-1}$ and $c>0$ is $A_r$-consistent.
Moreover, let $g=g_0+g_1$ denote any splitting of $g\in L^r_{0,\om}+L^r_{\om,1/\eta}$.
Defining the characteristic function $\chi_j$ of $\Si\cap\Si_j$ and the scalar
$$\begin{array}{rl}
m_j & = \di \frac{1}{|\Si\cap\Si_j|}\int_{\Si\cap\Si_j}
        (\varphi_jg_0+u'\cdot\na'\varphi_j+\be\vp_j u_n)dx'\\
    & = \di \frac{1}{|\Si\cap\Si_j|}\int_{\Si\cap\Si_j}(i\xi u_n-g_1)\varphi_jdx',
\end{array}$$
we split $g_j$ in the form
$$ g_j=g_{j0}+g_{j1}:=(\varphi_jg_0+u'\cdot \na'\varphi_j+\be\vp_j u_n-m_j\chi_j)+
   (\varphi_jg_1+m_j\chi_j).$$
Concerning $g_{j1}$ we get
$$ \begin{array}{l}
\|g_{j1}\|^r_{r,\om;\Si_j}  = \di\int_{\Si\cap\Si_j}
   |\varphi_jg_1+m_j|^r\om\, dx'\\[2ex]
 \hspace{1cm}\leq c(r)\big(\|g_1\|^r_{r,\om;\Si}+|m_j|^r\om(\Si\cap\Si_j)\big)\ek
\hspace{1cm} \leq c(r) \Big(\|g_1\|^r_{r,\om;\Si}+
   \di\frac{\om(\Si\cap\Si_j)\cdot\om'(\Si\cap\Si_j)^{r/{r'}}}{|\Si\cap\Si_j|^r}
   (\|\xi u_n\|^r_{(W^{1,r'}_{\om'})^*}+\|g_1\|^r_{r,\om;\Si})\Big)
\end{array}$$
with $c(r)>0$ independent of $\om$. Since we chose the balls
$B_j$ for $j\geq 1$ small enough, for each $j\geq 0$ there is a cube $Q_j$
with $\Si\cap\Si_j\subset Q_j$ and $|Q_j|<c(n)|\Si\cap\Si_j|$
where the constant $c(n)>0$ is independent of $j$. Therefore
\begin{equation}
\label{E4.5}
\begin{array}{rl}
\|g_{j1}\|_{r,\om;\Si_j}&\leq c(r)\Big(\|g_1\|_{r,\om} +
   \frac{c(n) \om(Q_j)^{1/r}\cdot\,\om'(Q_j)^{1/{r'}}}
    {|Q_j|}(\|\xi u_n\|_{(W^{1,r'}_{\om'})^*} +
   \|g_1\|_{r,\om})\Big)\ek
   & \leq c(r)(1+\car^{1/r})\big(\|\xi u_n\|_{(W^{1,r'}_{\om'})^*}+
     \|g_1\|_{r,\om;\Si}\big)
\end{array}
\end{equation}
for $j\geq 0$.
Furthermore, for every test function $\Psi\in C^{\infty}_0(\bar{\Si}_j)$ let
$$ \tilde{\Psi}=\Psi-\frac{1}{|\Si\cap\Si_j|}\int_{\Si\cap\Si_j}\Psi dx'.$$
By the definition of $m_j\chi_j$ we have $\int_{\Si_j}g_{j0}\,dx'=0$;
hence by Poincar\'e's inequality (see Proposition \ref{P2.7})
$$ \begin{array}{l}
|\int_{\Si_j}g_{j0}\Psi dx'|  =|\int_{\Si_j}g_{j0}\tilde{\Psi}dx'|\ek
 = |\int_{\Si}g_0(\varphi_j\tilde{\Psi})dx' +
    \int_{\Si}u'\cdot(\na'\varphi_j)\tilde{\Psi}dx'
    +\int_{\Si}\be u_n\varphi_j\tilde{\Psi}\,dx'|\ek
 \leq \|g_0\|_{-1,r,\om}\|\na'(\varphi_j\tilde{\Psi})\|_{r',\om'} +
        \|u'\|_{(W^{1,r'}_{\om'})^*}\|(\na'\varphi_j)\tilde{\Psi}\|_{1,r',\om'}
        +\|\be u_n\|_{(W^{1,r'}_{\om'})'}\|\vp_j \tilde{\Psi}\|_{1,r,\om'}\ek
 \leq c(\|g_0\|_{-1,r,\om}+\|u\|_{(W^{1,r'}_{\om'})^*})\|\na'\Psi\|_{r',\om';\Si_j},
\end{array}$$
where $c>0$ is $A_r$-consistent.  Thus
\begin{equation}
\label{E4.6}
\|g_{j0}\|_{-1,r,\om;\Si_j}\leq c\big(\|g_0\|_{-1,r,\om} +
\|u\|_{(W^{1,r'}_{\om'})^*}\big)\quad\text{for  } j\geq 0.
\end{equation}
Summarizing \eq{4.5} and \eq{4.6}, we get for $j\geq 0$
$$\|g_j;\hw^{-1,r}_\om(\Si_j)+L^r_{\om,1/\xi}(\Si_j)\|
\leq c\big(\|u\|_{(W^{1,r'}_{\om'})^*}+\|g; L^r_{0,\om}+ L^r_{\om,1/\xi}\|_0\big)
$$
with an $A_r$-consistent $c=c(r, \car)>0$, which yields in view of
$\xi^2>\be^2$ that
\begin{equation}
\label{E4.7}
\|g_j;\hw^{-1,r}_\om(\Si_j)+L^r_{\om,1/\xi}(\Si_j)\|
\leq c\big(\|u\|_{(W^{1,r'}_{\om'})^*}+\|g; L^r_{0,\om}+ L^r_{\om,1/\eta}\|_0\big)
\end{equation}
with an $A_r$-consistent $c=c(r, \car)>0$.

To complete the proof, apply Theorem \ref{T3.1} (i) to \eq{4.2}, \eq{4.3} when $j=0$.
Further use Theorem \ref{T3.1} (ii)
 in \eq{4.2}, \eq{4.3} for $j\geq 1$, but with $\la$ replaced by
$\la +M$ with $M=\la_0+\alpha_0$, where $\la_0=\la_0(\ve,r,\car)$
is the $A_r$-consistent constant indicated in Theorem \ref{T3.1} (ii).
This shift in $\la$ implies that $f_j$ has to be replaced by
$f_j+M\varphi_ju$ and that \eq{3.31} will be used with $\la$ replaced by $\la +M$.
Summarizing \eq{3.9}, \eq{3.31} as well as \eq{4.4}, \eq{4.7} and summing over all $j$
we arrive at \eq{4.1} with the additional terms
$$ I = \|Mu\|_{r,\om} + \|Mu\|_{(W^{1,r'}_{\om'})^*} +
   \|Mg; L^r_{0,\om}+ L^r_{\om,1/\eta}\|_0  $$
on the right-hand side of the inequality.
Note that $M=M(\ve,r,\car)$ is $A_r$-consistent, $|\eta|\leq \max\{\sqrt{2}|\xi|,\sqrt{2}\beta\}$
and that $g=\div'u'+i\eta u_n$ defines a natural splitting of
$g\in L^r_{0,\om}(\Si)+L^r_\om(\Si)$. Hence Poincar\'{e}'s inequality yields
$$ \begin{array}{rl}
I&\leq M\big(\|u\|_{r,\om;\Si}+\|\div'u'\|_{-1,r,\om}+\|u_n\|_{r,\om;\Si}\big)\ek
&\leq c_1\|u\|_{r,\om;\Si}\leq c_2\|\na' u\|_{r,\om;\Si}
\end{array}$$
with $A_r$-consistent constants $c_i=c_i(\ve,r,\Si,\car)>0,\, i=1,2$.

Thus \eq{4.1} is  proved.
\par\bigskip
Next, consider the case $\xi^2\leq \be^2$.
Since $S(u,p)=(f,-g)$, we have
\begin{equation}
\label{E3.3}
\begin{array}{l}
(\la-\Da')u'+\na' p=f'-\eta^2u', \,\,\div'u'=g-i\eta u_n,\quad\text{in }\Si,\ek
 u'|_{\pa\Si}=0,
\end{array}
\end{equation}
and
\begin{equation}
\label{E3.4}
\begin{array}{l}
(\la-\Da')u_n=f_n-\eta^2 u_n-i\eta p,\quad \text{in }\Si,\ek
 u_n|_{\pa\Si}=0.
\end{array}
\end{equation}
Now, apply \cite{Fr01}, Lemma 3.2 to
\eq{3.3}.   Then,  in view of $|\eta|\leq \sqrt{2}\be$
and  Poincar\'e's inequality, for all $\la\in -\al+S_\ve$, $\al\in (0,\al_0-\be^2)$ we have
$$\begin{array}{l}
\|(\la+\al)u', \na'^2u', \na'p\|_{r,\om;\Si}\ek
\leq c\big(\|f,\eta^2u\|_{r,\om;\Si}+|\la|\|g-i\eta u_n\|_{\hat{W}^{-1,r}_\om(\Si)}
+\|g-i\eta u_n\|_{{W}^{1,r}_\om(\Si)}+|\la|\|u'\|_{(W^{1,r'}_{\om'}(\Si))'}\big)\ek
\leq c\big(\|f, \na'u, p\|_{r,\om;\Si}+\|g\|_{{W}^{1,r}_\om(\Si)}+ |\la|\|g-i\eta u_n\|_{\hat{W}^{-1,r}_\om(\Si)}
+|\la|\|u\|_{(W^{1,r'}_{\om'}(\Si))'}\big)
\end{array}$$
with
$A_r$-consistent constant $c=c(r,\ve,\al, \be, \Si, {\cal A}_r(\Om))$.

In order to control $\|g-i\eta u_n\|_{\hat{W}^{-1,r}_\om(\Si)}$,
let us split $g$ as $g=g_0+g_1$, $g_0\in L^r_{0,\om}(\Si)$, $g_1\in L^r_{\om,1/\eta}(\Si)$.
Since $g_1-i\eta u_n$ has mean value zero in $\Si$, we get by poincar\'e's inequality
that
$$\begin{array}{l}
|\lan g_1-i\eta u_n, \psi\ran|=|\lan g_1-i\eta u_n, \bar\psi\ran|\ek
|\eta|\|g_1/\eta\|_{r,\om}\|\|\bar\psi\|_{r',\om'}
  +|\eta|\|u_n\|_{(W^{1,r'}_{\om'}(\Si))'}\|\bar\psi\|_{W^{1,r'}_{\om'}(\Si)}\ek
\leq c(r,\Si)(\|g_1/\eta\|_{r,\om}\|+\|u_n\|_{(W^{1,r'}_{\om'}(\Si))'})\|\na'\psi\|_{r',\om';\Si},
\end{array}
$$
for all $\psi\in C^\infty(\bar\Si)$, where $\bar\psi=\psi-\frac{1}{|\Si|}\int_{\Si} \psi\,dx'$.
Therefore,
$$
\|g-i\eta u_n\|_{\hat{W}^{-1,r}_\Om(\Si)}\leq \|g_0\|_{\hat{W}^{-1,r}_\Om(\Si)}
+c(\|g_1/\eta\|_{r,\om}\|+\|u_n\|_{(W^{1,r'}_{\om'}(\Si))'}).
$$
Thus, for all $\la\in -\al+S_\ve$, $\al\in (0,\al_0-\be^2)$ we have
\begin{equation}
\label{E3.5}
\begin{array}{l}
\|(\la+\al)u', \na'^2u', \na'p\|_{r,\om;\Si}\ek
\leq c\big(\|f, \na'u, p\|_{r,\om;\Si}+\|g\|_{{W}^{1,r}_\om(\Si)}+ |\la|\|u\|_{(W^{1,r'}_{\om'}(\Si))'}
+|\la|\|g: L^r_{0,\om}+L^r_{\om,1/\eta}\|_0\big)
\end{array}
\end{equation}
with
$A_r$-consistent constant $c=c(r,\ve,\al, \be, \Si, {\cal A}_r(\Om))$.

On the other hand, applying well-known results 
for the Laplace resolvent equations (cf. \cite{Fr01}) to \eq{3.4}, we get that
\begin{equation}
\label{E3.6}
\|(\la+\al)u_n, \na'^2u_n\|_{r,\om;\Si}
\leq c(\|f_n, u, p\|_{r,\om;\Si}
\end{equation}
 with $c=c(r,\ve,\al, \be, \Si, {\cal A}_r(\Om))$.
Thus, from \eq{3.5} and \eq{3.6} the assertion of the lemma for the case $\xi^2\leq \be^2$ is proved.

The proof of the lemma is complete. \qed

\vspace {0.1cm}
%
%
\begin{lem}
\label{L3.7}
Let $1<r<\infty,\, \om\in A_r$ and
$\xi\in\R^*$, $\be \in (0, \sqrt{\bar\al})$, $\al\in (0,\bar\al-\be^2)$,
 $\ve\in (0,\arctan{\frac{\sqrt{\bar\al-\be^2-\al}}{\be}})$,
$\la\in -\alpha+S_\ve$.
 Then there is an $A_r$-consistent constant
$c=c(\alpha,\ve,r,\be,\Si, \car)$ such that for every $(u,p)\in {\cal D}(S)$ and
$(f,-g)=S(u,p)$ the estimate
\begin{equation}
\label{E4.8}
\begin{array}{l}
\|\mu^2_+u,\mu_+\na'u,\na'^2u,\na'p,\eta p\|_{r,\om}\\[1ex]
\hspace{15mm} \leq c\big(\|f,\na'g,g,\xi g\|_{r,\om}+
(|\la|+1)\|g; L^r_{0,\om}+ L^r_{\om,1/\xi}\|_0\big)
\end{array}
\end{equation}
holds; here $\mu_+=|\la+\alpha+\xi^2|^{1/2}$.
\end{lem}
{\bf Proof:}
Assume that this lemma is wrong.  Then there  is a constant $c_0>0$, a sequence
$\{\om_j\}^\infty_{j=1}\subset A_r$ with ${\cal A}_r(\om_j)\leq c_0$ for all $j$,
sequences $ \{\la_j\}^\infty_{j=1}\subset
  -\alpha+S_\ve, \{\xi_j\}^\infty_{j=1}\subset \R^*$
and $(u_j,p_j)\in {\cal D}(S^{\om_j}_{r,\la_j,\xi_j})$ for all $j\in\N$ such that
\begin{equation}
\label{E4.9}
\begin{array}{l}
\|(\la_j+\alpha+\xi^2_j)u_j,(\la_j+\alpha+\xi^2_j)^{1/2}\na'u_j,
   \na'^2u_j,\na'p_j,\eta_jp_j\|_{r,\om_j}\ek
\hspace{2cm}\geq j\big(\|f_j,\na'g_j,g_j,\xi_jg_j\|_{r,\om_j} +
   (|\la_j|+1)\|g_j; L^r_{m,\om_j}+ L^r_{\om_j,1/\eta_j}\|_0
\end{array}
\end{equation}
where $\eta_j=\xi_j+i\be$, $(f_j,-g_j)=S^{\om_j}_{r,\la_j,\eta_j}(u_j,p_j)$.
Fix an arbitrary cube $Q$ containing $\Si$.
We may assume without loss of generality  that
\begin{equation}
\label{E4.10}
{\cal A}_r(\om_j)\leq c_0,\quad \om_j(Q)=1\quad \forall j\in \N,
\end{equation}
by using the $A_r$-weight $\tilde{\om}_j:=\om_j(Q)^{-1}\om_j$ instead of
$\om_j$ if necessary.
Note that \eq{4.10} also holds for $r',\{\om'_j\}$ in the following form:
${\cal A}_r(\om_j)\leq c_0^{r'/r},\; \om_j'(Q) \leq c_0^{r'/r} |Q|^{r'}.$
Therefore, by a minor modification of Proposition \ref{P2.5} (3),
there exist numbers $s, s_1$ such that
\begin{equation}
\label{E4.11}
L^r_{\om_j}(\Si)\hookrightarrow L^s(\Si),\quad L^{s_1}(\Si)
\hookrightarrow L^{r'}_{{\om}'_j}(\Si),\quad j\in \N,
\end{equation}
with embedding constants independent of $j\in\N$.
Furthermore, we may assume without loss of generality that
\begin{equation}
\label{E4.12}
\|(\la_j+\alpha+\xi^2_j)u_j,(\la_j+\alpha+\xi^2_j)^{1/2}
   \na'u_j,\na'^2u_j,\na'p_j,\eta_jp_j\|_{r,\om_j}=1
\end{equation}
and consequently that
\begin{equation}
\label{E4.13}
\|f_j,\na'g_j,g_j,\xi_jg_j\|_{r,\om_j} +
(|\la_j|+1)\|g_j; L^r_{m,\om_j}+ L^r_{\om_j,1/\xi_j}\|_0\rightarrow 0
    \quad \text{as} \quad j\rightarrow \infty.
\end{equation}
From \eq{4.11}, \eq{4.12} we have
\begin{equation}
\label{E4.14}
\|(\la_j+\alpha+\xi^2_j)u_j,(\la_j+\alpha+\xi^2_j)^{1/2}
   \na'u_j,\na'^2u_j,\na'p_j,\eta_jp_j\|_s \leq K,
\end{equation}
with some $K>0$ for all $j\in\N$ and
\begin{equation}
\label{E4.15}
\|f_j,\na'g_j,g_j,\eta_jg_j\|_s\rightarrow 0 \quad\text{as}
   \quad j\rightarrow \infty.
\end{equation}
Without loss of generality let us suppose that as $j\rightarrow \infty$,
$$\begin{array}{l}
\la_j\rightarrow \la\in -\alpha+\bar{S}_\ve \quad \text{or}
   \quad|\la_j|\rightarrow \infty\\
\xi_j\rightarrow 0  \quad \text{or}\quad \xi_j\rightarrow \xi \neq 0
   \quad \text{or}\quad|\xi_j|\rightarrow \infty.
\end{array}$$
Thus we have to consider \vspace{0.2cm} six possibilities.

{\it (i) The case} $\la_j\rightarrow \la \in -\alpha+\bar{S}_\ve,\quad
\xi_j\rightarrow \xi < \infty$. \\
Due to \eq{4.14} $\{u_j\}\subset W^{2,s}$ and $\{p_j\}\subset W^{1,s}$
are bounded sequences. In virtue of the compactness of the embedding
$W^{1,s}(\Si)\hookrightarrow L^s(\Si)$ for the bounded domain $\Si$,
we may assume (suppressing indices for subsequences) that
\begin{equation}
\label{E4.16}
\begin{array}{lcl}
u_j\rightarrow u, \na'u_j\rightarrow \na'u \quad & \text{in }L^s
& \quad \text{(strong convergence)}\\[1ex]
\na'^2u_j \wrar \na'^2u \quad &\text{in } L^s & \quad \text{(weak convergence)}\\[1ex]
p_j\rightarrow p \quad &\text{  in  }L^s &\quad \text{(strong convergence)}\\[1ex]
 \na'p_j \wrar \na'p \quad &\text{in }  L^s & \quad \text{(weak convergence)}
\end{array}
\end{equation}
for some $(u,p)\in {\cal D}(S_{s,\la,\xi})$ as $j\rightarrow \infty$.
Therefore, $S_{s,\la,\xi}(u,p)=0$ and, consequently, $u=0,\,p=0$ by Lemma \ref{L3.4}.
On the other hand we get from \eq{4.12} that
$\sup_{j\in\N}\|u_j\|_{2,r,\om_j}<\infty$   and
$\sup_{j\in\N}\|p_j\|_{1,r,\om_j}<\infty$
which, together with the weak convergences $u_j\wrar 0$ in $W^{2,s}(\Si)$,
$p_j\wrar 0$ in  $W^{1,s}(\Si)$, yields
$$ \|u_j\|_{1,r,\om_j}\rightarrow 0,\quad \|p_j\|_{r,\om_j}\rightarrow 0$$
due to Proposition \ref{P2.6} (2).
Moreover, since $\sup_{j\in\N}\|\la_j u_j\|_{r,\om_j}<\infty$
and $\la_j u_j\wrar \la u=0$ in $L^s(\Si)$, Proposition \ref{P2.6} (3) implies that
\begin{equation} \label{E4.16n}
\|\la_j u_j\|_{(W^{1,r'}_{\om'_j})^*}\rightarrow 0.
\end{equation}
Thus \eq{4.1}, \eq{4.12} and \eq{4.13} yield the contradiction $1\leq 0$.

{\it (ii) The case} $\la_j \ra \la \in -\alpha+\bar{S_\ve},
\quad |\xi_j|\rightarrow \infty$.\\
 From \eq{4.12} we get
$\|\na' u_j,\xi_j u_j,p_j\|_{r,\om_j}\rightarrow 0.$
On the other hand, since   $\|u_j\|_{r,\om_j} \rightarrow 0$ and
$u_j\rightarrow 0$ in $L^s$ as $j\rightarrow \infty$, Proposition \ref{P2.6} (3)
implies \eq{4.16n}. Thus, from \eq{4.1}, \eq{4.12} and \eq{4.13}
we get the contradiction $1\leq 0$.
\par\bigskip
{\it (iii) The case} $|\la_j|\rightarrow \infty,\quad \xi_j\rightarrow \xi<\infty$.\\
By \eq{4.12}
\begin{equation}
\label{E4.22}
\|\na'u_j,\xi_ju_j\|_{r,\om_j}\rightarrow 0\quad \text{as}\quad j\rightarrow \infty.
\end{equation}
Further, \eq{4.14} yields the convergence
$$\begin{array}{lcl}
u_j\rightarrow 0, \na'u_j\rightarrow 0 & \quad \text{and}\quad
& \na'^2u_j\wrar 0, \la_ju_j\wrar v,\\
 p_j\rightarrow p&\quad\text{and} \quad & \na'p_j\wrar \na'p,
\end{array}$$
in $L^s$, which, together with \eq{4.15}, leads to
\begin{equation}
\label{E4.23}
v'+\na'p=0,\quad v_n+i\eta p=0.
\end{equation}

Let $g_j:=g_{j0}+g_{j1}$, $g_{j0}\in L^r_{0,\om_j}$,
$g_{j1}\in L^r_{\om_j}$. Then, by \eq{4.14} we have
\begin{equation}
\label{E4.17}
\|\la_j g_{j0}\|_{-1,r,\om_j}+\|\la_j g_{j1}/\eta_j\|_{r,\om_j}\ra 0\quad (j\ra\infty).
\end{equation}
From \eq{4.11}, \eq{4.17} we see that
\begin{equation} \begin{array}{rl}
|\lan\la_jg_j,\varphi\ran|
& = |\lan\la_jg_{j0},\varphi\ran + \lan\la_jg_{j1},\varphi\ran| \nn\ek
& \leq \|\la_jg_{j0}\|_{-1,r,\om_j}\|\na'\varphi\|_{r',\om'_j} +
  \|\la_jg_{j1}\|_{r,\om_j}\|\varphi\|_{r',\om'_j}\label{E4.23n} \ek
& \leq c \big(\|\la_jg_{j0}\|_{-1,r,\om_j}\| +
  \|\la_jg_{j1}/\eta_j\|_{r,\om_j}\big)\|\varphi\|_{W^{1,s_1}(\Si)}. \nn
\end{array} \end{equation}
Consequently,
\begin{equation}
\label{E4.24}
\la_jg_j\in(W^{1,s_1}(\Si))^*\quad\text{and}
\quad\| \la_jg_j\|_{(W^{1,s_1}(\Si))^*}\rightarrow 0
\quad \text{as }j\rightarrow \infty.
\end{equation}
Therefore, it follows from the divergence equation $\div'_{\eta_j}u_{j} = g_{j}$
that for all $\varphi\in C^\infty(\bar{\Si})$
$$ \begin{array}{rcl}
\lan v',-\na'\varphi\ran + \lan i\eta v_n,\varphi \ran
& = & \lim_{j\rightarrow \infty}\lan\div'\la_j u'_j +
      i\la_j\xi_ju_{jn},\varphi\ran\ek
& = & \lim_{j\rightarrow \infty}\lan\la_jg_j,\varphi\ran=0 ,
\end{array}$$
yielding $\div'v' = -i\eta v_n,\; v'\cdot N|_{\pa\Si}=0.$
Therefore \eq{4.23} implies
\begin{equation}
\label{E3.7}
-\Da'p+\eta^2p=0\text{ in }\Si, \quad \di \frac{\pa p}{\pa N}=0
   \text{ on }\pa\Si.
\end{equation}
Here note that $\eta^2=\xi^2-\be^2+2i\xi\be$.
Hence, if $\xi\neq 0$ then $p\equiv 0$ since the all eigenvalues of the
Neumann Laplacian in $\Si$ is real; if $\xi=0$, then
$\eta^2=-\be^2$ and hance $p\equiv 0$ due to the condition
$\be^2<\bar\al\leq \al_1$. That is,  we have $p\equiv 0$, and
 and also $v\equiv 0$.
Now, due to Proposition \ref{P2.6} (2), (3), we get \eq{4.16n} and the convergence
$\|p_j\|_{r,\om_j}\rightarrow 0,$
since $\la_ju_j\wrar 0$ in $L^s$, $p_j\wrar0$ in $W^{1,s}$ and
$\sup_{j\in\N} \|\la_ju_j\|_{r,\om_j}<\infty,$
$\sup_{j\in\N} \|p_j\|_{1,r,\om_j}<\infty$.
Thus \eq{4.1}, \eq{4.12}, \eq{4.13} and \eq{4.22} lead to the contradiction
$1\leq 0$.

{\it (iv) The case} $|\la_j|\rightarrow \infty, \quad |\xi_j|\rightarrow \infty.$\\
To come to a contradiction, it is enough to prove \eq{4.16n} since
$\|\na'u_j,\xi_ju_j,p_j\|_{r,\om_j}\rightarrow 0$ as $j\rightarrow \infty$.
From \eq{4.12} we get the convergence
$$ \begin{array}{lcl}
   u_j\rightarrow 0, \na'u_j\rightarrow 0 & \quad\text{and}\quad
&  \na'^2u_j\wrar 0, (\la_j+\eta^2_j)u_j\wrar v,\\
   p_j\rightarrow 0 & \quad\text{and}\quad
&  \na'p_j\wrar 0, \quad  \eta_jp_j\wrar q
\end{array}$$
in $L^s$ with some $v,q\in L^s$.
Therefore, \eq{4.15} and $(R_{\la_j,\xi_j})$ yield
$$ v'=0, \quad v_n+iq=0.$$
Since $\|\la_ju_j\|_s\leq c_\ve \|(\la_j+\eta^2_j)u_j\|_s$, there exists
$w=(w',w_n)\in L^s$ such that, for a suitable subsequence, $\la_ju_j\wrar w$.
Let $g_j=g_{j0}+g_{j1},\,j\in\N,$ be a sequence of splittings satisfying \eq{4.17}.
By \eq{4.11} we get for all $\varphi\in C^\infty(\bar{\Si})$
$$ |\lan\la_jg_{j0},\varphi\ran|  +
   \Big|\lan\di\frac{\la_jg_{j1}}{\eta_j},\varphi\ran\Big|
   \ra 0  \quad \text{as}\quad j\rightarrow \infty,$$
cf. \eq{4.23n} and \eq{4.24}.
Hence, the divergence equation implies that for $j\ra\infty$
$$ \lan\la_ju_{jn},\varphi\ran\ek
   =  \di\frac{1}{i\eta_j}\lan\la_jg_{j0},\varphi\ran
   + \lan\frac{\la_jg_{j1}}{i\eta_j},\varphi\ran
   + \frac{1}{i\eta_j}\lan\la_ju'_j,\na'\varphi\ran \ra 0$$
for all $\varphi \in C^{\infty}(\bar{\Si})$ yielding $\lan w_n,\varphi\ran =0$ and
consequently $w_n=0$.

Obviously, $\eta_ju_j\rightarrow 0$ in $L^s$ as $j\rightarrow \infty$.
Therefore, by \eq{4.15} and the boundedness of the sequence
$\big\{ \|\eta_j \na u_j\|_{r,\om_j} \big\}$, we get from the identity
$\div'(\eta_ju'_j) + i\eta^2_ju_{jn} = \eta_jg_j$ that
$$\eta^2_ju_{jn}\wrar 0 \,\,\text{ and hance }\,\xi^2_ju_{jn}\wrar 0\quad\text { in } L^s \text{ as } j\rightarrow \infty.$$
Thus we proved $v_n=0$. Now $v=0$ together with the estimate
$\|(\la_j+\xi^2_j)u_j\|_{r,\om_j}\leq 1$ imply due to Proposition \ref{P2.6} (3) that
$ \|(\la_j+\xi^2_j)u_j\|\rightarrow 0$ in $(W^{1,r'}_{\om'_j})^*$
as $j\rightarrow \infty.$ Hence also \eq{4.16n} is proved.  \vspace{0.2cm}

Now the proof of this lemma is complete. \hfill \qed \vspace{0.4cm}
%
%
\begin{theo}
\label{T3.8}
Let $1<r<\infty,\, \om\in A_r$ and
$\xi\in\R^*$, $\be \in (0, \sqrt{\bar\al})$, $\al\in (0,\bar\al-\be^2)$,
 $\ve\in (0,\arctan{\frac{\sqrt{\bar\al-\be^2-\al}}{\be}})$.
Then for every $\la\in-\alpha+S_\ve$, $\xi\in \R^*$ and
$f\in L^r_\om(\Si),\, g\in W^{1,r}_\om(\Si)$ the parametrized resolvent problem
$(R_{\la,\xi,\be})$ has a unique solution
$(u,p)\in \big(W^{2,r}_\om(\Si)\cap W^{1,r}_{0,\om}(\Si)\big)\ti W^{1,r}_\om(\Si)$.
Moreover, this solution satisfies the estimate \eq{4.8}
with an $A_r$-consistent constant $c=c(\alpha,\be,\ve,r,\Si,\car)>0$.
\end{theo}
{\bf Proof:} The existence is obvious since, for every
$\la\in-\alpha+S_\ve,\xi\in \R^*$ and $\om\in A_r(\R^{n-1})$, the range
${\cal R}(S^\om_{r,\la,\xi})$ is closed and dense in
$L^r_\om(\Si)\ti W^{1,r}_\om(\Si)$ by Lemma \ref{L3.5} and by Lemma \ref{L3.4}, respectively.
Here note that for fixed $\la\in \C,\,\xi\in\R^*$ the norm
$\|\na'g,g,\xi g\|_{r,\om}+ (1+|\la|) \|g; L^r_{m,\om}+ L^r_{\om,{1/\xi}}\|_0$
is equivalent to the norm of $W^{1,r}_\om(\Si)$.
The uniqueness of solutions is obvious from Lemma \ref{L3.4}. \vspace{0.4cm}\hfill\qed

Now, for fixed $\om\in A_r, 1<r<\infty$, define the
operator-valued functions
$$\begin{array}{l}
a_1: \R^*\rightarrow {\cal L}(L^r_\om(\Si);W^{2,r}_{0,\om}(\Si)
     \cap W^{1,r}_\om(\Si)),\\[1ex]
b_1: \R^*\rightarrow  {\cal L}(L^r_\om(\Si);W^{1,r}_\om(\Si))
\end{array}$$
by
\begin{equation}
\label{E4.26}
 a_1(\xi)f:=u_1(\xi),\quad b_1(\xi)f:=p_1(\xi),
\end{equation}
where $(u_1(\xi), p_1(\xi))$ is the solution to $(R_{\la,\xi,\be})$
corresponding to $f\in L^r_\om(\Si)$ and $g=0$. Further, define
$$\begin{array}{l}
a_2: \R^*\rightarrow {\cal L}(W^{1,r}_\om(\Si);W^{2,r}_{0,\om}(\Si)
     \cap W^{1,r}_\om(\Si)),\\[1ex]
b_2: \R^*\rightarrow {\cal L}(W^{1,r}_\om(\Si);W^{1,r}_\om(\Si))
\end{array}$$
by
\begin{equation}
\label{E4.27}
a_2(\xi)g:=u_2(\xi),\quad b_2(\xi)g:=p_2(\xi).
\end{equation}
with $(u_2(\xi), p_2(\xi))$ the solution to $(R_{\la,\xi,\be})$ corresponding to
$f=0$ and $g\in W^{1,r}_\om(\Si)$.
%
\begin{coro}
\label{C3.9}
Assume the same for $\alpha,\be,\xi,\la$ as in Theorem \ref{T3.8}.
Then, the operator-valued functions $a_1, b_1$ and $a_2, b_2$  defined by
\eq{4.26}, \eq{4.27} are Fr\'{e}chet differentiable in $\xi\in\R^*$.
Furthermore, their derivatives
$w_1 = \frac{d}{d\xi}a_1(\xi)f,\; q_1 = \frac{d}{d\xi}b_1(\xi)f$ for fixed
$f\in L^r_\om(\Si)$ and
$w_2=\frac{d}{d\xi}a_2(\xi)g,\; q_2=\frac{d}{d\xi}b_2(\xi)g$ for fixed
$g\in W^{1,r}_\om(\Si)$ satisfy the estimates
\begin{equation}
\label{E4.28}
\|(\la+\alpha)\xi w_1,\xi\na'^2w_1, \xi^3w_1, \xi\na'q_1,\xi\eta q_1\|_{r,\om}
\leq c\|f\|_{r,\om}
\end{equation}
and
\begin{equation}
 \label{E4.29}
\begin{array}{l}
\|(\la+\alpha)\xi w_2,\xi\na'^2w_2, \xi^3w_2, \xi\na'q_2,\xi\eta q_2\|_{r,\om}\\[1ex]
\hspace{15mm} \leq c
\big(\|\na'g,g,\xi g\|_{r,\om}+(|\la|+1)\|g; L^r_{0,\om}+ L^r_{\om,1/\eta}\|_0\big),
\end{array}
\end{equation}
with an $A_r$-consistent constant $c=c(\alpha, \be, r,\ve,\Si,\car)$ independent of
$\la\in -\alpha+S_\ve$ and $\xi\in\R^*$.
\end{coro}
{\bf Proof:} Since $\xi$ enters in $(R_{\la,\xi})$ in a polynomial way,
it is easy to prove that $a_j(\xi),b_j(\xi), j=1,2,$ are Fr\'{e}chet
differentiable and their derivatives $w_j,q_j$ solve the system
\begin{equation}
\label{E4.30}
\begin{array}{rcl}
(\la+\eta^2-\Da')w'_j+\na'q_j&=&-2\eta u'_j\\[1ex]
(\la+\eta^2-\Da')w_{jn}+i\eta q_j&=&-2\eta u_{jn}-ip_j\\[1ex]
\div'w'_j+i\eta w_{jn}&=&-iu_{jn},
\end{array}
\end{equation}
where $(u_1,p_1), (u_2,p_2)$ are the solutions to $(R_{\la,\xi,\be})$ for
$f\in L^r_\om(\Si),\, g=0$ and $f=0,\, g\in W^{1,r}_\om(\Si)$, respectively.

We get from \eq{4.30} and Theorem \ref{T3.8} for $j=1,2,$
\begin{equation}
\label{E4.31}
\begin{array}{l}
\|(\la+\alpha)\xi w_j,\xi\na'^2w_j, \xi^3w_j, \xi\na'q_j,\xi\eta q_j\|_{r,\om}\\[1ex]
\hspace{1cm} \leq c\big(\|\xi\eta u'_j,\xi p_j, \xi\na'u_{jn}, \xi^2 u_{jn}\|_{r,\om}
 + (|\la|+1)\|i\eta u_{jn}; L^r_{0,\om}+ L^r_{\om,1/\eta}\|_0\big)\\[1ex]
\hspace{1cm} \leq c\big(\|\xi^2 u_j,\xi p_j, \xi\na' u_j\|_{r,\om}
 + (|\la|+1)\|u_j\|_{r,\om}\big)\ek
\hspace{1cm} \leq c \|u_j, (\la+\alpha+\xi^2)u_j,
   \sqrt{\la+\alpha+\xi^2}\na' u_j, \xi p_j\|_{r,\om}\ek
\hspace{1cm} \leq c \|(\la+\alpha+\xi^2)u_j,
   \sqrt{\la+\alpha+\xi^2}\na' u_j, \na'^2 u_j,\xi p_j\|_{r,\om},
\end{array}
\end{equation}
with an $A_r$-consistent constant $c=c(\alpha,r,\ve,\Si,\car)$;
here we used the fact that
$\xi^2+ |\la+\alpha|\leq c(\ve,\alpha) |\la+\alpha+\xi^2|$  for all
$\la\in -\alpha+S_\ve, \xi\in \R$ and
$\|u_j\|_{r,\om}\leq c(\car)\|\na'^2u_j\|_{r,\om}$
(see \cite{Fr01}, Corollary 2.2). Thus Theorem \ref{T3.8} and \eq{4.31}
prove \eq{4.28}, \eq{4.29}.  \hfill \qed

\section{Proof of the Main Results}
\subsection{Proof of Theorem \ref{T2.1} -- Theorem \ref{T2.3} }
\noindent
The proof of Theorem 2.1 is based on the theory of operator-valued
Fourier multipliers.
The classical  H\"{o}rmander-Michlin theorem for scalar-valued multipliers
for $L^q(\R^k),\, q\in(1,\infty),\, k\in\N,$ extends to an operator-valued
version for Bochner spaces $L^q(\R^k;X)$ provided that $X$ is a {\it UMD space}
and that the boundedness condition for the derivatives of the multipliers
is strengthened to {\it ${\cal R}$-boundedness}.
%
%
\begin{tdefi}
\label{D4.1}
A Banach space $X$ is called a {\it UMD} space
if the Hilbert transform
$$ Hf(t)=-\frac{1}{\pi}\,\,\text{PV}\,\, \int\frac{f(s)}{t-s}\,ds
\quad\text{ for  }f\in{\cal S}(\R;X),$$
where ${\cal S}(\R;X)$ is the Schwartz space of all rapidly
decreasing $X$-valued functions, extends to a bounded linear operator
in $L^q(\R;X)$ for some $q \in(1,\infty)$.
\end{tdefi}

It is well known that, if $X$ is a {\it UMD} space, then the
Hilbert transform is bounded in $L^q(\R;X)$ {\it for all}
$q\in(1,\infty)$ (see e.g. \cite{RRT86}, Theorem 1.3) and that
weighted Lebesgue spaces $L^r_\om(\Si), 1<r<\infty,\,\om\in A_r,$
are {\it UMD} spaces.
%
%
\begin{tdefi}
\label{D4.2}
Let $X,Y$ be Banach spaces. An operator family ${\cal T}\subset{\cal L}(X;Y)$
is called $\cal R$-bounded if there is a constant $c>0$ such that for all
$T_1,\ldots,T_N\in{\cal T},$ $x_1,\ldots, x_N\in X$ and $N\in\N$
\begin{equation}
\label{E5.1}
\big\|\sum_{j=1}^{N}\ve_j(s)T_jx_j\big\|_{L^q(0,1;Y)}
\leq c\, \big\|\sum_{j=1}^{N}\ve_j(s)x_j\big\|_{L^q(0,1;X)}
\end{equation}
for some $q\in [1,\infty)$, where $(\ve_j)$ is any sequence of independent,
symmetric $\{-1,1\}$-valued random variables on $[0,1]$.
The smallest constant $c$ for which {\rm \eq{5.1}} holds is denoted by
$R_q({\cal T})$, the $\cal R$-bound of ${\cal T}$.
\end{tdefi}

\begin{rem}
\label{R4.3}
{\rm
(1) Due to {\it Kahane's inequality} (\cite{DJT95})
\begin{equation}
\label{E5.2}
\big\|\sum_{j=1}^{N}\ve_j(s)x_j\big\|_{L^{q_1}(0,1;X)}
\leq c(q_1,q_2,X) \big\|\sum_{j=1}^{N}\ve_j(s)x_j\big\|_{L^{q_2}(0,1;X)},\,\,\,
1\leq q_1,q_2<\infty,
\end{equation}
the inequality \eq{5.1} holds {\it for all } $q\in [1,\infty)$
if it holds for some $q\in [1,\infty)$.

(2) If an operator family ${\cal T}\subset{\cal L}(L^r_\om(\Si))$,
$1<r<\infty,\;\om\in A_r(\R^{n-1})$, is $\cal R$-bounded, then
${\cal R}_{q_1}({\cal T})\leq C {\cal R}_{q_2}({\cal T})$ for all
$q_1,q_2\in [1,\infty)$ with a constant $C=C(q_1,q_2)>0$
independent of $\om$. In fact, introducing the isometric isomorphism
$$ I_\om: L^r_\om(\Si)\ra  L^r(\Si),\quad I_\om f=f\om^{1/r},$$
for all $T\in {\cal L}(L^r_\om(\Si))$ we have
$\tilde{T}_\om=I_\om T I_\om^{-1}\in {\cal L}(L^r(\Si))$ and
$\|T\|_{{\cal L}(L^r_\om(\Si))}=\|\tilde{T}_\om\|_{{\cal L}(L^r(\Si))}$.
Then it is easily seen that ${\cal \tilde{T}}_\om
:=\{I_\om T I_\om^{-1}: \,\, T\in {\cal T}\} \subset{\cal L}(L^r(\Si))$
is ${\cal R}$-bounded and ${\cal R}_q({\cal \tilde{T}}_\om)={\cal R}_q({\cal T})$
for all $q\in [1,\infty)$. Thus the assertion follows.
}
\end{rem}

\begin{tdefi}
\label{D4.4}
(1) Let $X$ be a Banach space and $(x_n)^\infty_{n=1}\subset X$.
A series $\sum_{n=1}^{\infty}x_n$ is called unconditionally convergent
if $\sum_{n=1}^{\infty}x_{\si(n)}$ is convergent in norm for every
permutation $\si: \N\rightarrow \N$.

(2) A sequence of projections $(\Da_j)_{j\in\N}\subset {\cal L}(X)$ is
called a Schauder decomposition of a Banach space $X$ if
$$ \Da_i\Da_j=0\;\; \text{for all}\;\; i\neq j,\quad
   \sum_{j=1}^{\infty}\Da_jx=x \;\;\text{for each} \;\; x\in X.$$
A Schauder decomposition $(\Da_j)_{j\in\N}$ is called unconditional
if the series $\sum_{j=1}^{\infty}\Da_jx$ converges unconditionally
for each $x\in X$.
\end{tdefi}
%
\begin{rem}
\label{R4.5}
{\rm
(1) If $(\Da_j)_{j\in\N}$ is an unconditional Schauder decomposition of a
Banach space $Y$, then for each $p\in [1,\infty)$ there is a constant
$c_\Da=c_\Da(p)>0$ such that for all $x_j$ in the range
${\cal R}(\Da_j)$ of $\Da_j$ the inequalities
\begin{equation}
\label{E5.8}
c^{-1}_\Da \Big\|\sum_{j=l}^{k}x_j\Big\|_Y\quad \leq
\quad \Big\|\sum_{j=l}^{k}\ve_j(s) x_j\Big\|_{L^p(0,1;Y)}
\leq c_\Da  \Big\|\sum_{j=l}^{k}x_j\Big\|_Y
\end{equation}
are valid for any sequence $(\ve_j(s))$ of independent, symmetric
$\{-1,1\}$-valued random variables defined on $(0,1)$ and
for all $l\leq k\in\Z$, see e.g. \cite{DHP03}, (3.8).

(2) Let $Y=L^q(\R;L^r_\om(\Si))$ and assume that each $\Da_j$
commutes with the isomorphism $I_\om$ introduced in Remark \ref{R4.3} (2).
Then the constant $c_\Da$ is easily seen to be
independent of the weight $\om$.

(3) In the previous definitions and results the set of indices
$\N$ may be replaced by $\Z$ without any further changes.

(4) Let $X$ be a {\it UMD} space and $\chi_{[a,b)}$ denote the characteristic
function for the interval $[a,b)$. Let $R_s=\finv\chi_{[s,\infty)}\fou$ and
$$\Da_j:=R_{2^j}-R_{2^{j+1}},\; j \in \Z.$$
It is well known that the Riesz projection $R_0$ is bounded in
$L^q(\R;X)$ and that the set $\{R_s-R_t: \,\, s,t\in \R\}$ is
$\cal R$-bounded in ${\cal L}(L^q(\R;X))$ for each $q\in(1,\infty)$.
In particular, $\{\Da_j:\; j\in\Z\}$ is  $\cal R$-bounded in ${\cal L}(L^q(\R;X))$
and an unconditional Schauder decomposition of $R_0 L^q(\R;X)$, the image of
$L^q(\R;X)$ by the Riesz projection $R_0$, see \cite{DHP03},
proof of Theorem 3.19.
}
\end{rem}

We recall an operator-valued Fourier multiplier theorem
in Banach spaces.
Let ${\cal D}_0(\R;X)$ denote the set of $C^\infty$-functions
$f: \R\rightarrow X$ with compact support in $\R^*$.

%
\begin{theo} {\em (\cite{DHP03}, Theorem 3.19, \cite{We01}, Theorem 3.4)}
\label{T4.6}
Let $X$ and $Y$ be UMD spaces and $1<q<\infty$.
Let $M: \R^*\rightarrow {\cal L}(X,Y)$ be a differentiable function
such that
$${\cal R}_q \big(\{M(t),\,tM'(t):\,\,t\in\R^*\}\big) \leq A.$$
Then the operator
$$Tf=\big(M(\cdot)\hat{f}(\cdot)\big)^\vee,\quad f \in {\cal D}_0(X),$$
 extends to a bounded operator $T: L^q(\R;X)\rightarrow L^q(\R;Y)$
with operator norm $\|T\|_{{\cal L}(L^q(\R;X);L^q(\R;Y))}\leq CA$
where $C>0$ depends only on $q,X$ and $Y$.
\end{theo}

\begin{rem}
\label{R4.7}
Checking the proof of \cite{DHP03}, Theorem 3.19, one can see
that the constant $C$ in Theorem \ref{T4.6} equals
$$C={\cal R}({\cal P})\cdot (c_\Da)^2$$
where ${\cal R}({\cal P})$ is the $\cal R$-bound of the operator family
${\cal P}=\{R_s-R_t: \,\,s,t \in \R\}$  in ${\cal L}(L^q(\R;{X}))$
and $c_\Da$ is the {\it unconditional constant} of the
{\it Schauder decomposition} $\{\Delta_j:\,\, j\in\Z\}$
of the space $R_0L^q(\R;{X})$; see \cite{DHP03}, Section 3, for details.
In particular, for $X=L^r_\om(\Si),\, 1<r<\infty,\,\om\in A_r,$ using
the isometry $I_\om$ of Remark \ref{R4.3} (2), we get that the constants
${\cal R}({\cal P})$, see Remark \ref{R4.3} (2), and $c_\Da$ do not depend
on the weight $\om$; concerning $c_\Da$ we again use that $I_\om$
commutes with each $\Delta_j.$
\end{rem}
       %
\begin{theo} {\em (Extrapolation Theorem)}
\label{T4.8}
Let $1<r,s<\infty,\, \om\in A_r(\R^{n-1})$ and $\Si\subset \R^{n-1}$ be
an open set. Moreover let ${\cal T}$ be
a family of linear operators with the property that there exists an
$A_s$-consistent constant $C_{\cal T}=C_{\cal T}({\cal A}_s(\nu))>0$
such that for all $\nu\in A_s$
$$\|Tf\|_{s,\nu}\leq C_{\cal T}\|f\|_{s,\nu}$$
for all $T\in{\cal T}$ and all $f\in L^s_\nu(\Si)$.
Then every $T\in {\cal T}$ can be extended to $L^r_\om(\Si)$ and
${\cal T}$ is $\cal R$-bounded in ${\cal L}(L^r_\om(\Si))$ with an
$A_r$-consistent $\cal R$-bound $c_{\cal T}(q,r,\car)$, i.e.,
\begin{equation}
\label{E5.5}
{\cal R}_q({\cal T})\leq c_{\cal T}(q,r,\car)\quad
   \text{for all}\quad q\in (1,\infty).
\end{equation}
\end{theo}
{\bf Proof:}  From the proof of \cite{Fr01}, Theorem 4.3, it can be deduced
that ${\cal T}$ is ${\cal R}$-bounded in ${\cal L}(L^r_\om(\Si))$ and that \eq{5.5}
is satisfied for $q=r$. Then, Remark \ref{R4.3} yields \eq{5.5} for every
$1<q<\infty$.\qed
\par\bigskip
Now we are in a position to prove Theorem 2.1.\vspace{0.3cm}\\
%
{\bf Proof of Theorem 2.1:}
 Let $f(x',x_n):= e^{\be x_n}F(x',x_n)$ for
 $(x',x_n)\in\Si\ti\R$
 and let us define $u,p$ in the cylinder
$\Om=\Si\ti \R$ by
$$ u(x)={\cal F}^{-1}(a_1 \hat{f})(x),
   \quad  p(x)={\cal F}^{-1}(b_1 \hat{f})(x),$$
where $a_1, b_1$ are the operator-valued multiplier functions defined
in \eq{4.26}.
We will show that  $(U, P)=(e^{\be x_n}u, e^{\be x_n} p)$ is the unique solution to $(R_\la)$
with $g=0$ such that
\begin{equation}
\label{E5.6}
(u, p) \in \big(W^{2;q,r}_\om(\Om)\cap W^{1;q,r}_{0,\om}(\Om)\big)
\ti \hw^{1;q,r}_\om (\Om)
\end{equation}
and the estimate \eq{2.1} holds. Obviously, $(U,P)$ solves the resolvent problem
$(R_\la)$ with $g=0$.
For $\xi\in \R^*$ define $m_\la(\xi):L^r_\om(\Si)\rightarrow L^r_\om(\Si)$ by
$$ m_\la(\xi)f := \big((\la+\alpha)a_1(\xi)\hat{f}, \xi\na'a_1(\xi)\hat{f},
   {\na'}^2 a_1(\xi)\hat{f}, \xi ^2a_1(\xi)\hat{f},
   \na'b_1(\xi)\hat{f},\xi b_1(\xi)\hat{f}\big).$$
Theorem \ref{T3.8} and Corollary \ref{C3.9} show that the operator family
$\{m_\la(\xi), \xi m'_\la(\xi):\xi\in \R^*\}$ satisfies the assumptions of
Theorem \ref{T4.8}, e.g., with $s=r$. Therefore, this operator family is
$\cal R$-bounded in ${\cal L}(L^r_\om(\Si))$; to be more precise,
$$ {\cal R}_q\big(\{m_\la(\xi), \xi m'_\la(\xi):\xi\in \R^*\}\big)
   \leq c(q,r,\alpha,\be, \ve,\Si,\car) < \infty.$$
Hence Theorem \ref{T4.6} and Remark \ref{R4.7} imply that
$$ \|(m_\la\hat{f})^\vee\|_{L^q(L^r_\om)} \leq C\|f\|_{L^q(L^r_\om)}$$
with an $A_r$-consistent constant $C=C(q,r,\alpha,\be, \ve,\Si,\car)>0$
independent of the resolvent parameter $\la\in -\alpha+S_\ve$.
Note that, due to the definition of the multiplier $m_\la(\xi)$,
we have $(\la+\alpha)u,\na^2u, \na p\in L^q(L^r_\om)$ and
$$ \|(\la+\alpha)u,\na^2 u,\na p\|_{L^q(L^r_\om)}
   \leq \|(m_\la\hat{f})^\vee\|_{L^q(L^r_\om)}.$$
Thus the existence of a solution satisfying \eq{2.1} is proved.

The uniqueness of solutions is obvious by the uniqueness result
for $\be=0$ of \cite{FR05-3}, Theorem 2.1.
 Now the proof of Theorem 2.1
is complete. \hfill \qed\\[2ex]
%
%
 {\bf Proof of Corollary 2.2:} Defining the Stokes operator $A=A_{q,r;\be,\om}$
 by \eq{2.2}, due to the Helmholtz decomposition of the space
 $L^q_\be(L^r_\om)$ on the cylinder $\Om$, see \cite{Fa03},
 we get that for $F\in L^q_\be(L^r_\om)_\si$ the solvability of the equation
\begin{equation}
\label{E5.7}
(\la+A)U=F \quad \text{in}\quad L^q_\be(L^r_\om)_\si
\end{equation}
is equivalent to the solvability of $(R_\la)$ with right-hand side
$G\equiv 0$. By virtue of Theorem 2.1 for every
$\la\in -\alpha+S_\ve$ there exists a unique solution
$U=(\la+A)^{-1}F\in D(A)$ to \eq{5.7} satisfying the estimate
$$ \|(\la+\alpha)U\|_{L^q_\be(L^r_\om)_\si}
=\|(\la+\alpha)u\|_{L^q(L^r_\om)}
\leq C\|f\|_{L^q(L^r_\om)}=C\|F\|_{L^q_\be(L^r_\om)_\si}$$
with $C=C(q,r,\alpha,\be,\ve,\Si,\car)$ independent of $\la$,
where $u=e^{\be x_n}U$, $f=e^{\be x_n}F$.
Hence \eq{2.3} is proved.
Then \eq{2.4} is a direct consequence of \eq{2.3} using semigroup theory.
\qed
\par\bigskip
\noindent
{\bf Proof of Theorem 2.3:} The proof will be  done if we show that the operator family
$${\cal T}=\{\la(\la+A_{q,r;\be,\om})^{-1}: \;\la\in \mathrm{i}\R\}$$
is $\cal R$-bounded in ${\cal L}(L^q_\be(L^r_\om)_\si)$.
By the way, since $L^q_\be(L^r_\om)_\si$ is isomorphic to a closed
subspace $X$ of $L^q(L^r_\om)$
with isomorphism $I_\be F:= e^{\be x_n}F$,
it is enough to show ${\cal R}$-boundedness of
$${\cal \tilde{T}}=\{I_\be\la(\la+A_{q,r;\be,\om})^{-1}I_\be^{-1}: \;\la\in \mathrm{i}\R\}
\subset {\cal L}(X).$$
In the following we write shortly
$$H_\be\equiv I_\be\la(\la+A_{q,r;\be,\om})^{-1}I_\be^{-1}.$$

For $\xi\in\R^*$ and $\la\in S_\ve$,
 let $m_\la(\xi):=\la a_1(\xi)$ where $a_1(\xi)$
is the solution operator for $(R_{\la,\xi,\be})$ with $g=0$ defined by \eq{4.26}.
Then, we have
 $$H_\be f= I_\be\la U=\la I_\be U=(m_\la(\xi)\hat{f})^\vee,
                  \quad \forall f\in{\cal S}(\R; L^r_\om(\Si))\cap X,$$
with $U$ the solution to $(R_\la)$ with $F=I_\be^{-1}f$, $G=0$.
Note that
${\cal S}(\R; L^r_\om(\Si))$ is dense in $L^q(\R; L^r_\om(\Si))$
and hence ${\cal S}(\R; L^r_\om(\Si))\cap X$ is dense in $X$.
Hence, in view of Definition \ref{D4.2} and Remark \ref{R4.3},
${\cal R}$-boundedness of
${\cal \tilde{T}}$ in ${\cal L}(X)$
is proved if there is a constant $C>0$ such that
\begin{equation}
\label{E5.9}
\big\|\sum_{i=1}^N\ve_i(m_{\la_i}\hat{f}_i)^\vee
      \big\|_{L^q(0,1;L^q(\R: L^r_\om(\Si)))}
\leq C \big\|\sum_{i=1}^N\ve_i f_i\big\|_{L^q(0,1;L^q(\R: L^r_\om(\Si)))}
\end{equation}
for any independent, symmetric and $\{-1,1\}$-valued random variables
$(\ve_i(s))$ defined on $(0,1)$, for all $(\la_i)\subset \mathrm{i} \R$
and $(f_i)\subset {\cal S}(\R; L^r_\om(\Si))\cap X$.
Without loss of generality we may assume that
$\supp \hat{f}_i\subset [0,\infty)$, $i=1,\ldots,N$,
since $R_0f:=(\chi_{[0,\infty)}(\xi)\hat{f})^\vee$ is continuous in $L^q(\R; L^r_\om(\Si))$ and
 $$f_i(x',x_n)=(\chi_{[0,\infty)}\hat{f}_i(\xi))^\vee(x',x_n)+
(\chi_{[0,\infty)}\hat{f}_i(-\xi))^\vee(x',-x_n).$$
Note that, if $\supp \hat{f}\subset [0,\infty)$,
then $\supp(m_\la \hat{f}) \subset [0,\infty)$ as well.
Therefore, instead of \eq{5.9}
we shall prove the estimate
\begin{equation}
\label{E5.9n}
\big\|\sum_{i=1}^N\ve_i(m_{\la_i}\hat{f}_i)^\vee
      \big\|_{L^q(0,1;Y)}
\leq C \big\|\sum_{i=1}^N\ve_i f_i\big\|_{L^q(0,1;Y)}
\end{equation}
 for $(f_i)\subset {\cal S}(\R; L^r_\om(\Si))\cap X\cap Y$.

Obviously $m_\la(\xi) = m_\la(2^j)+\int_{2^j}^{\xi}m'_\la(\tau)\,d\tau$
for $\xi\in [2^j,2^{j+1}),\, j\in \Z,$ and
$\big(m_\la(2^j)\widehat{\Da_jf}\big)^\vee{}=m_\la(2^j)\Da_jf$ for
$f\in {\cal S}(\R; L^r_\om(\Si))\cap X\cap Y$. Furthermore,
$$\begin{array}{l}
\di\Big(\int_{2^j}^{\xi}m'_\la(\tau)\,d\tau\, \widehat{\Da_jf}(\xi)\Big)^\vee  =
        \Big(\di\int_{2^j}^{2^{j+1}}m'_\la(\tau)
        \chi_{[2^j,\xi)}(\tau)\widehat{\Da_jf}(\xi)\,d\tau\Big)^\vee \\[1ex]
\hspace{2cm} = \di \Big(\int_{0}^{1}2^jm'_\la(2^j(1+t)) \chi_{[2^j,\xi)}(2^j(1+t))
        \chi_{[2^j,2^{j+1})}(\xi)\hat{f}(\xi)\,dt\Big)^\vee \\[1ex]
\hspace{2cm} = \di \int_{0}^{1}2^jm'_\la(2^j(1+t))B_{j,t}\Da_jf\, dt.
\end{array}$$
where $B_{j,t}=R_{2^j(1+t)}-R_{2^{j+1}}$. Thus we get
\begin{equation}
\label{E5.10}
\begin{array}{rcl}
\di\big(m_\la(\xi)\hat{f}(\xi)\big)^\vee & = &
\di\sum_{j\in\Z}\Big(\big(m_\la(2^j)+\int_{2^j}^{\xi}m'_\la(\tau)
   \,d\tau\big)\,\widehat{\Da_jf}\Big)^\vee\ek
 &  = &\di\sum_{j\in\Z} \big(m_\la(2^j)\widehat{\Da_jf}\big)^\vee{}+
                \sum_{j\in\Z} \Big(\int_{2^j}^{\xi}m'_\la(\tau)
   \,d\tau\,\widehat{\Da_jf}\Big)^\vee\ek
& = &\di\sum_{j\in\Z} m_\la(2^j)\Da_jf +
                \sum_{j\in\Z} \int_{0}^{1}2^jm'_\la(2^j(1+t))B_{j,t}\Da_jf\, dt.
\end{array}
\end{equation}

First let us prove
\begin{equation}
\label{E5.11}
\big\|\sum_{i=1}^N\ve_i(s) \sum_{j\in \Z}m_{\la_i}(2^j)\Da_jf_i\big\|_{L^q(0,1;Y)}
\leq C \big\|\sum_{i=1}^N \ve_i(s) f_i\big\|_{L^q(0,1;Y)}.
\end{equation}
Note that the operator $m_{\la_i}(2^j)$ commutes with $\Da_j$, $j\in\Z$;
hence, for almost all $s\in (0,1)$, the sum
$\sum_{i=1}^{N}\ve_i(s)m_{\la_i}(2^j)\Da_jf_i$
belongs to the range of $\Da_j$.
Therefore, for any $l,k\in\Z$ we get by \eq{5.8} that
\begin{equation}
\label{E5.12}
\begin{array}{l}
\di\big\|\sum_{i=1}^N\ve_i
\sum_{j=l}^{k}m_{\la_i}(2^j)\Da_jf_i\big\|_{L^q(0,1;Y)} \ek
\hspace{2cm} = \di\Big(\int_0^1\big\|\sum_{j=l}^k
  \sum_{i=1}^{N}\ve_i(s)m_{\la_i}(2^j)\Da_jf_i\big\|^q_Y\,ds\Big)^{1/q} \ek
\hspace{2cm} \leq c_\Da\di \Big(\int_0^1\int_0^1\big\|\sum_{j=l}^k \ve_j(\tau)
\sum_{i=1}^{N}\ve_i(s)m_{\la_i}(2^j)\Da_jf_i\big\|^q_Y\,d\tau\,ds\Big)^{1/q}\ek
\hspace{2cm} = c_\Da\di \big\|\sum_{i=1}^{N}
\sum_{j=l}^k \ve_{ij}(s,\tau) m_{\la_i}(2^j)\Da_jf_i\big\|_{L^q((0,1)^2;Y)}
\end{array}
\end{equation}
where $\ve_{ij}(s,\tau)=\ve_i(s)\ve_j(\tau)$;
note that $(\ve_{ij})_{i,j\in\Z}$ is a sequence of
independent, symmetric and $\{-1,1\}$-valued random variables
defined on $(0,1)\ti(0,1)$. Furthermore, due to Theorem \ref{T3.8}, the operator family
$\{m_\la(\xi):\,\la\in \mathrm{i}\R,\,\xi\in\R^*\}
\subset {\cal L}(L^r_\om(\Si))$
is uniformly bounded by an $A_r$-consistent constant,
and hence it is $\cal R$-bounded by Theorem \ref{T4.8}. Therefore, using
Fubini's theorem and \eq{5.8}, we proceed in \eq{5.12} as follows:
\begin{equation}
\label{E5.12n}\begin{array}{l}
=c_\Da\di \big\|\sum_{i=1}^{N}\sum_{j=l}^k \ve_{ij}(s,\tau)
        m_{\la_i}(2^j)\Da_j f_i\big\|_{L^q(\R;L^q((0,1)^2;L^r_\om(\Si)))}\ek
\leq Cc_\Da\di \big\|\sum_{i=1}^{N}\sum_{j=l}^k
     \ve_{ij}(s,\tau) \Da_j f_i\big\|_{L^q(\R;L^q((0,1)^2;L^r_\om(\Si)))}\ek
= Cc_\Da\di \big\|\sum_{i=1}^{N}\sum_{j=l}^k
  \ve_{ij}(s,\tau) \Da_j f_i\big\|_{L^q((0,1)^2;Y)}
\leq Cc^2_\Da\di \big\|\sum_{i=1}^{N} \ve_i \sum_{j=l}^k
     \Da_j f_i\big\|_{L^q(0,1;Y)}.
\end{array}
\end{equation}
Since $\{\sum_{j=l}^k\Da_j:\,\, l,k\in\Z\}$ is $\cal R$-bounded in ${\cal L}(Y)$ and
$(\Da_j)$ is a Schauder decomposition of $Y$, we see by Lebesgue's theorem that
the right-hand side of \eq{5.12n} converges to $0$ as either
$l,k\ra \infty$ or $l,k\ra -\infty$.
Thus, by \eq{5.12}, \eq{5.12n}, the series
$\sum_{i=1}^N\ve_i(s) \sum_{j\in\Z}m_{\la_i}(2^j)\Da_j f_i$
converges in $L^q(0,1;Y)$, and \eq{5.11} holds.

Next let us show that
\begin{equation}
\label{E5.13}
\big\|\sum_{i=1}^N\ve_i(s) \sum_{j\in\Z}\int_0^1
        2^jm'_{\la_i}(2^j(1+t))B_{j,t}\Da_j f_i\, dt\big\|_{L^q(0,1;Y)}
\leq C \big\|\sum_{i=1}^N\ve_i(s) f_i\big\|_{L^q(0,1;Y)}.
\end{equation}
Using the same argument as in the proof of \eq{5.11} and the
$\cal R$-boundedness of the operator families
$\{B_{j,t}:\,\, j\in\Z, t\in (0,1)\}\subset {\cal L}(Y)$ and
$ \{2^j(1+t)m'_{\la}(2^j(1+t)):\,\,\la\in \mathrm{i}\R, j\in \Z,
  t\in (0,1)\}\subset {\cal L}(L^r_\om(\Si))$,
see Corollary \ref{C3.9}, we have
$$\begin{array}{l}
\di \big\|\sum_{i=1}^N\ve_i(s) \sum_{j=l}^k\int_0^1
     2^jm'_{\la_i}(2^j(1+t))B_{j,t}\Da_j f_i\,dt\big\|_{L^q(0,1;Y)}\ek
\hspace{2cm} \leq \di\int_0^1 \big\|\sum_{i=1}^N\ve_i(s)
     \sum_{j=l}^k 2^jm'_{\la_i}(2^j(1+t))B_{j,t}\Da_j f_i\big\|_{L^q(0,1;Y)}\,dt\ek
\hspace{2cm} \leq \di c_\Da\int_0^1 \big\|\sum_{i=1}^N\sum_{j=l}^k \ve_{ij}(s,\tau)
        2^jm'_{\la_i}(2^j(1+t))B_{j,t}\Da_j f_i\big\|_{L^q((0,1)^2;Y)}\,dt\ek
\hspace{2cm} \leq \di c_\Da\int_0^1 \big\|\sum_{i=1}^N\sum_{j=l}^k \ve_{ij}(s,\tau)
        2^j(1+t)m'_{\la_i}(2^j(1+t))\Da_j f_i\big\|_{L^q((0,1)^2;Y)}\,dt\ek
\hspace{2cm} \leq Cc_\Da^2 \di \big\|\sum_{i=1}^{N} \ve_i(s)
     \sum_{j=l}^k\Da_jf_i\big\|_{L^q((0,1);Y)}
\end{array}$$
for all $l,k\in\Z.$ Thus \eq{5.13} is proved.

By \eq{5.11}, \eq{5.13} we conclude that the operator family
${\cal T}=\{\la(\la+A_{q,r;\be,\om})^{-1}: \;\la\in \mathrm{i}\R\}$
is $\cal R$-bounded in ${\cal L}(L^q_\be(L^r_\om))$.
Then, by \cite{We01}, Corollary 4.4, for each
$f\in L^p(\R_+; L^q_\be(L^r_\om)_\si), 1<p<\infty,$
the mild solution $U$ to the system
\begin{equation}
\label{E5.14}
U_t+A_{q,r;\be, \om}U=F,\quad u(0)=0
\end{equation}
belongs to
$L^p(\R_+; L^q_\be(L^r_\om)_\si)\cap L^p(\R_+; D(A_{q,r;\be,\om}))$
and satisfies the estimate
$$ \| U_t, A_{q,r;\be,\om}U\|_{L^p(\R_+; L^q_\be(L^r_\om)_\si)}
   \leq C\|F\|_{L^p(\R_+; L^q_\be(L^r_\om)_\si)}.$$
Furthermore, \eq{2.3} with $\la=0$ implies that $U$ also satisfies
this inequality. If $F\in L^p(\R_+; L^q_\be(L^r_\om))$, let
$U$ be the solution of \eq{5.14} with $F$ replaced by $P_{q,r;\be,\om}F$, where
$P_{q,r;\be,\om}$ denotes the Helmholtz projection in
$L^q_\be(L^r_\om))$, and define $P$ by $\na P=(I-P_{q,r;\be,\om})(f-u_t+\Delta u)$.
By \eq{2.1} with $\la=0$ and the boundedness of $P_{q,r;\be,\om}$ we get \eq{2.5c}.
Finally, assume $e^{\al t} F\in L^p(\R_+; L^q_\be(L^r_\om))$
for some $\al\in(0,\bar\al-\be^2)$ and let $V$ be the solution of the system
$V_t+(A-\al)V = e^{\al t}P_{q,r;\be,\om} F,\; V(0)=0$.
Obviously, replacing $A$ by $A-\al$ in the previous arguments,
$v$ is easily seen to satisfy estimate \eq{2.5b}. Then
$U(t)=e^{-\al t}V(t)$ solves \eq{5.14} and satisfies \eq{2.5d}.
In each case the constant $C$ depends only on $\car$ due to Remark \ref{R4.7}.

The proof of Theorem 2.3 is complete. \hfill\qed
\par\bigskip
\subsection{Proof of Theorem \ref{T2.4} and Theorem \ref{T2.5}}
{\bf Proof of Theorem \ref{T2.4}:}
Let $1<q<\infty$ and
$\xi\in\R^*$, $\be \in (0, \sqrt{\al^*})$, $\al^*=\min_{1\leq i\leq m}{\bar\al_i}$,
$\al\in (0,\al^*-\be^2)$,
 $\ve\in (0,\arctan{\frac{\sqrt{\al^*-\be^2-\al}}{\be}})$.
Fix $\la\in-\alpha+S_\ve$ and $\xi\in \R^*$.
Note that $\la+A_{q,\tb}$ with $\be_i=0$ for all $i=1,\ldots,m$
is injective and surjective, see \cite{FR05-4}, Theorem 1.2.
Hence,
given any $F\in  L^q_{\tb,\si}(\Om)$, for all $\la\in -\al+S_\ve$
there is a unique
$(U, \na P)\in D(A_q)\ti L^q(\Om)$ such that
\begin{equation}
\label{E4.40}
 \begin{array}{rcll}
\la U- \Da U +  \na P & =& F & \mbox{ in } \Om,\ek
\div U & = & 0 & \mbox{ in }\Om,\ek
U & = & 0 & \mbox{ on } \pa \Om,
\end{array}
\end{equation}

Without loss of generality we may assume that there exist cut-off functions
$\{\varphi_i\}_{i=0}^{m}$ such that
\begin{equation}
\label{E4.15n}
\begin{array}{l}
    \sum_{i=0}^{m}\varphi_i(x)=1, \quad 0\leq \varphi_i(x)\leq 1\quad
    \text{for }x\in\Om,\ek
\varphi_i \in C^{\infty}(\bar{\Om}_i), \quad
   \text{dist}\,(\text{supp}\,\varphi_i,\, \pa\Om_i\cap \Om)\geq \delta>0,
   \,\,i=0,\ldots, m,
\end{array}
\end{equation}
where 'dist' means the distance. In what follows, for $i=1,\ldots, m$ let $\widetilde{\Om}_i$ be
the infinite straight cylinder extending the semi-infinite cylinder $\Om_i$, and
denote the zero extension of $\vp_i v$ to $\widetilde{\Om}_i$ by $\widetilde{\vp_iv}$.

Let
$$(u^{0},p^{0}): =(\vp_0 U,\vp_0 P),\;
(u^{i}, p^{i}):=(\widetilde{\vp_i U}, \widetilde{\vp_i P})\;\text{ for }i=1,\ldots,m.$$
Then $(u^{0},p^{0})$ on $\Om_0$ satisfies
$$ \begin{array}{rcll}
\la u^0- \Da u^0 +  \na p^{0} & =& f^0 & \mbox{ in } \Om_0,\ek
\div u^0 & = & g^0 & \mbox{ in }\Om_0,\ek
u^0 & = & 0 & \mbox{ on } \pa \Om_0,
\end{array}$$
and $(u^i, p^i)$ on
$\widetilde{\Om}_i,\,i=1,\ldots, m$, satisfy
$$ \begin{array}{rcll}
\la u^i- \Da u^i + \na p^i & =&\tilde{f}^i & \mbox{ in } \widetilde{\Om}_i,\ek
\div u^i & = & \tilde{g}^i & \mbox{ in }\widetilde{\Om}_i,\ek
u^i & = & 0 & \mbox{ on } \pa \widetilde{\Om}_i,
\end{array}$$
where
$$ f^i:=\varphi_iF+(\na\varphi_i)P-(\Da\varphi_i)U-2\na\varphi_i\cdot\na U,
   \quad g^i:=\na\varphi_i\cdot U,\quad i=0,\ldots, m.$$
Since  $\,\,\text{supp}\, g^i\subset \Om_0$, $g^i\in W^{1,q}_0(\Om_0)$
and $\int_{\Om_0}g^i\,dx=0$ for  $i=0,\ldots,m$, we get by the well-known theory of the
divergence problem for $i=0,\ldots,m$ that
there is some $w_i\in W^{2,q}_0(\Om_0)$ satisfying
$\div w_i=g^i$ in $\Om_0$
 and
\begin{equation}
\begin{array}{l}
\|\na^2 w_i\|_{L^q(\Om_0))}\leq c\|\na g^i\|_{L^q(\Om_0)}
\leq c\|\na U\|_{L^{q}_0(\Om_0)},\ek
\|w_i\|_{L^q(\Om_0))}\leq c\|g^i\|_{(W^{1,q}(\Om_0))^*}
\leq c\|U\|_{(W^{1,q}(\Om_0))^*},\ek

\end{array}
\end{equation}
where $c=c(\Om_0,q)$, cf. \cite{Ga94-1}.
Then $\tilde{w}_i$, the extension by $0$ of $w_i$ to $\widetilde{\Om}_i$,
$i=1\ldots,m,$ satisfies
\begin{equation}
\label{E4.17n}
\begin{array}{l}
e^{\be_i x^i_n}\na^2 \tilde{w}_i
\in L^q(\tilde\Om_i),\;
 \| e^{\be_i x^i_n}\na^2 \tilde{w}_i\|_{L^q(\tilde\Om_i)}
 \leq c\|\na U\|_{L^{q}(\Om_0)}.
\end{array}
\end{equation}

Now, $v^0:=u^0-w_0$ and $v^i:=u^i-\tilde{w}_i$, $i=1,\ldots,m$, solve, respectively,
$$ \begin{array}{rcll}
\la v^0- \Da v^0 + \na p^0 & =&f^0
-(\la w_0-\Da w_0) & \mbox{ in } \Om_0,\ek
\div v^0 & = & 0 & \mbox{ in }\Om_0,\ek
v^0 & = & 0 & \mbox{ on } \pa \Om_0,
\end{array}$$
and
$$ \begin{array}{rcll}
\la v^i- \Da v^i + \na p^i & =&\tilde{f}^i
-(\la\tilde{w}_i - \Da \tilde{w}_i) & \mbox{ in } \widetilde{\Om}_i,\ek
\div v^i & = & 0 & \mbox{ in }\widetilde{\Om}_i,\ek
v^i & = & 0 & \mbox{ on } \pa \widetilde{\Om}_i.
\end{array}$$
Then, using the fact that the Stokes operator in $L^q$-spaces on bounded
domains is injective and surjective we get that
\begin{equation}
\begin{array}{l}
\|v^0, \la v^0,\na^2 v^0, \na p^0\|_{L^q(\Om_0)}
\leq
 c\|F,\na U,P\|_{L^{q}(\Om_0)}+(|\la|+1)\|U\|_{(W^{1,q}(\Om_0))^*}
 \end{array}
\end{equation}
with $c$ independent of $\la$.
Moreover, by Theorem \ref{T2.1} we have
\begin{equation}
\begin{array}{l}
\|v^i, \la v^i, \na^2 v^i, \na p^i\|_{L^q_{\be_i}(\R; L^q(\Si^i))}
\leq
 c(\|F\|_{L^{q}_{\be_i}(\tilde\Om_i)}
+\|\na U, P\|_{L^{q}(\Om_0)}
\ek\hspace{2cm}+(|\la|+1)\|U\|_{(W^{1,q}(\Om_0))^*}),\; i=1,\ldots,m,
\end{array}
\end{equation}
with $c$ independent of $\la$.
Due to $U=\sum_{i=0}^m u^i$, $P=\sum_{i=0}^m p^i$ in $\Om$ and \eq{4.17n}, we get
$\na^2 U, \na P\in L^q_{\tb}(\Om)$ and
\begin{equation}
\label{E4.34}
\begin{array}{l}
 \|U, \la U,\na^2 U, \na P\|_{L^q_{\tb}(\Om)}\leq
  c(\|F\|_{L^{q}_{\tb}(\Om)}+\|\na U, P\|_{L^{q}(\Om_0)})
  \ek\hspace{5cm} +(|\la|+1)\|U\|_{(W^{1,q}(\Om_0))^*}.
\end{array}
\end{equation}
Indeed, by a contradiction argument \eq{4.34} yields
\begin{equation}
\label{E4.35}
\begin{array}{l}
 \|U, \la U,\na^2 U, \na P\|_{L^q_{\tb}(\Om)}\leq
  c(\|F\|_{L^{q}_{\tb}(\Om)}+\|\na U, P\|_{L^{q}(\Om_0)})
\end{array}
\end{equation}
with $c$ independent of $\la$.
\par Assume that \eq{4.35} does not hold.
Then there are sequences $\{\la_j\}\subset -\alpha+S_\ve,$
$\{(U_j,P_j)\}\in $ such that
\begin{equation}
\label{E4.37}
\|U_j, \la_j U_j,\na^2 U_j,\na P_j\|_{L^q_\tb(\Om)}=1,
    \quad\|F_j\|_{L^q_\tb(\Om)}\ra 0 \quad\text{as  }j\ra\infty,
\end{equation}
where $F_j=\la U_j-\Da U_j+\na P_j$.
Without loss of generality we may assume that
\begin{equation}
\label{E4.38}
\la_j U_j\wrar V,\; U_j\wrar U, \; \na^2 u_j\wrar \na^2 U, \;
\na P_j\wrar \na P\quad\text{as  }j\ra\infty
\end{equation}
with some $V\in L^q_\tb(\Om),\, U\in W^{2,q}_\tb(\Om)\cap W^{1,q}_{0,\tb}(\Om)$
 and $P\in \hw^{1,q}_\tb(\Om)$.
Moreover, we may assume $\int _{\Om_0} P_j \, dx=0,\, \int _{\Om_0} P \, dx=0$
and that $\la_j\ra \la\in \{-\alpha+\bar{S}_\ve\}\cup\{\infty\}$. \vspace{0.2cm}

\par (i) Let $\la_j\ra \la \in -\alpha+\bar{S}_\ve$.

Then, $V=\la U$ and
it follows that  $(U,P)$ solves \eq{4.40} with $F=0$
yielding $(U,P)=0$.
On the other hand, we have the strong convergence
\begin{equation}
\label{E4.39}
U_j\ra 0 \,\,\,\text{in  } W^{1,q}(\Om_0),\,\,\, P_j\ra 0\,\,\,\text{in  }L^q(\Om_0),\,\,\,
  (|\la_j|+1)U_j \ra 0\,\,\,\text{in  }(W^{1,q'}(\Om_0))^*
\end{equation}
due to the compact embeddings
$W^{2,q}(\Om_0)\subset\subset W^{1,q}(\Om_0)\subset\subset L^q(\Om_0)
                     \subset\subset (W^{1,q'}(\Om_0))^*$,
Poincar\'e's inequality on $\Om_0$.
Thus \eq{4.35} yields the contradiction $1\leq 0$.

\par (ii) Let $|\la_j|\ra \infty$.
Then, we conclude that
$U=0$, and consequently $V+\na P=0$ where $V\in L^q_{\si}(\Om)$.
Note that this is the $L^q$-Helmholtz decomposition of the null vector field on $\Om$.
Therefore, $V=0,\, \na P=0$.
Again we get \eq{4.39} and finally the contradiction $1\leq 0$.

The proof of the theorem is complete.
\qed
\par\bigskip\noindent
{\bf Proof of Theorem \ref{T2.5}:}
The idea of the proof is also to use a cut-off technique.
Note that any $F\in L^p(\R_+; L^q_\tb(\Om))$
also belongs to $L^p(\R_+; L^q(\Om))$ for $1<p,q<\infty$.
Hence, by maximal $L^p$-regularity
of the Stokes operator in $L^q(\Om)$, which follows by \cite{FR05-4}, Theorem 1.2,
we get that the problem \eq{2.7} has a unique solution $(U,\na P)$ such that
$$(U,\na P)\in
L^p(\R_+; W^{2,q}(\Om)\cap W^{1,q}_{0}(\Om)\cap L^q_\si(\Om))
\ti L^p(\R_+; L^q(\Om)), U_t\in  L^p(\R_+; L^q(\Om)).$$
We shall prove that this solution $(U, \na P)$, furthermore, satisfies
\begin{equation}
\label{E4.33}
(U,\na P)\in
L^p(\R_+; W^{2,q}_\tb(\Om))
\ti L^p(\R_+; L^q_\tb(\Om)), U_t\in L^p(\R_+; L^q_\tb(\Om)).
\end{equation}
Once \eq{4.33} is proved, the (linear) solution operator
$$L^p(\R_+; L^q_\tb(\Om))\ni F\mapsto
(U,\na P)\in L^p(\R_+; W^{2,q}_\tb(\Om)\cap W^{1,q}_0(\Om)\cap L^q_\si(\Om))
\ti L^p(\R_+; L^q_\tb(\Om))$$
is obviously closed and hence bounded by the closed graph theorem
thus implying \eq{2.6}.

The proof of \eq{4.33} is based on cut-off technique using Theorem \ref{T2.3}.
Let $\{\varphi_i\}_{i=0}^{m}$ be cut-off functions given by \eq{4.15n} and
let
$$(u^{0},p^{0}): =(\vp_0 U,\vp_0 P),\;
(u^{i}, p^{i}):=(\widetilde{\vp_i U}, \widetilde{\vp_i P})\;\text{ for }i=1,\ldots,m.$$
Then $(u^{0},p^{0})$ on $\Om_0$ satisfies
$$ \begin{array}{rcll}
u^{0}_t - \Da u^0 +  \na p^{0} & =& f^0 & \mbox{ in } \R_+\ti\Om_0\ek
\div u^0 & = & g^0 & \mbox{ in }\R_+\ti\Om_0\ek
u^0(0,x) & = & 0 & \mbox{ in } \Om_0,\ek
u^0 & = & 0 & \mbox{ on } \pa \Om_0,
\end{array}$$
and $(u^i, p^i)$ on
$\widetilde{\Om}_i,\,i=1,\ldots, m$, satisfy
$$ \begin{array}{rcll}
u^i_t - \Da u^i + \na p^i & =&\tilde{f}^i & \mbox{ in } \R_+\ti\widetilde{\Om}_i\ek
\div u^i & = & \tilde{g}^i & \mbox{ in }\R_+\ti\widetilde{\Om}_i\ek
u^i(0,x) & = & 0 & \mbox{ in } \widetilde{\Om}_i,\ek
u^i & = & 0 & \mbox{ on } \pa \widetilde{\Om}_i,
\end{array}$$
where
$$ f^i:=\varphi_iF+(\na\varphi_i)P-(\Da\varphi_i)U-2\na\varphi_i\cdot\na U,
   \quad g^i:=\na\varphi_i\cdot U,\quad i=0,\ldots, m.$$
Note that one has  $\,\,\text{supp}\, g^i\subset \Om_0$ hence $g^i\in L^p(\R_+; W^{1,q}_0(\Om_0))$
and $\int_{\Om_0}g^i\,dx=0$ for  $i=0,\ldots,m$.
Therefore, by the well-known theory of the
divergence problem for $i=0,\ldots,m$ there is some $w_i\in L^p(\R_+; W^{2,q}_0(\Om_0))$ such that
$\div w_i(t)=g^i(t)$ in $\Om_0$ for almost all $t\in\R_+$,
$w_{it}\in L^p(\R_+; L^{q}(\Om_0))$
 and
\begin{equation}
\label{E4.18}
\begin{array}{l}
\|\na^2 w_i\|_{L^p(\R_+; L^q(\Om_0))}\leq c\|\na g^i\|_{L^p(\R_+; L^q(\Om_0))}
\leq c\|\na U\|_{L^p(\R_+; L^{q}_0(\Om_0))},\ek
\|w_{it}\|_{L^p(\R_+; L^q(\Om_0))}\leq c\|g^i_t\|_{L^p(\R_+; (W^{1,q'}(\Om_0))')}
\leq c\|U_t\|_{L^p(\R_+; (W^{1,q'}(\Om_0))')},
\end{array}
\end{equation}
where $c=c(\Om_0,q)$, cf. \cite{Ga94-1}.
Then $\tilde{w}_i$, the extension by $0$ of $w_i$ to $\widetilde{\Om}_i$,
$i=1\ldots,m,$ satisfies
\begin{equation}
\label{E4.19}
\begin{array}{l}
e^{\be_i x^i_n}\tilde{w}_{it}, e^{\be_i x^i_n}\na^2 \tilde{w}_i
\in L^p(\R_+; L^q(\tilde\Om_i)),\ek
 \|e^{\be_i x^i_n}\tilde{w}_{it},
 e^{\be_i x^i_n}\na^2 \tilde{w}_i\|_{L^p(\R_+; L^q(\tilde\Om_i))}
 \leq c(\|\na U\|_{L^p(\R_+; L^{q}(\Om_0))}+\|U_t\|_{L^p(\R_+; (W^{1,q'}(\Om_0))')}).
\end{array}
\end{equation}
Moreover, note that $w_i(0,x)=0$ due to $g^i(0,x)=0$ for $x\in \Om$.

Now, $v^0:=u^0-w_0$ and $v^i:=u^i-\tilde{w}_i$, $i=1,\ldots,m$, solve, respectively,
$$ \begin{array}{rccl}
v^0_t - \Da v^0 + \na p^0 & =&f^0-w_{0t}
+\Da w_0 & \mbox{ in } \R_+\ti\Om_0,\ek
\div v^0 & = & 0 & \mbox{ in }\R_+\ti\Om_0,\ek
v^0(0,x) & = & 0 & \mbox{ in } \Om_0,\ek
v^0 & = & 0 & \mbox{ on } \pa \Om_0,
\end{array}$$
and
$$ \begin{array}{rccl}
v^i_t - \Da v^i + \na p^i & =&\tilde{f}^i-\tilde{w}_{it}
+\Da \tilde{w}_i & \mbox{ in } \R_+\ti\widetilde{\Om}_i,\ek
\div v^i & = & 0 & \mbox{ in }\R_+\ti\widetilde{\Om}_i,\ek
v^i(0,x) & = & 0 & \mbox{ in } \widetilde{\Om}_i,\ek
v^i & = & 0 & \mbox{ on } \pa \widetilde{\Om}_i.
\end{array}$$
Then, by the maximal regularity of Stokes operator in bounded
domains  in view of \eq{4.18} we obtain that
\begin{equation}
\label{E4.20}
\begin{array}{l}
\|v^{0},v^{0}_{t}, \na^2 v^0, \na p^0\|_{L^p(\R_+; L^q(\Om_0))}
\leq
 c(\|F,\na U,P\|_{L^p(\R_+; L^{q}(\Om_0))}
 +\|U_t\|_{L^p(\R_+; (W^{1,q'}(\Om_0))')}),
\end{array}
\end{equation}
and, by Theorem \ref{T2.3} in view of \eq{4.19}, that
\begin{equation}
\label{E4.21}
\begin{array}{l}
\|v^{i},v^{i}_t, \na^2 v^i, \na p^i\|_{L^p(\R_+; L^q_{\be_i}(\R; L^q(\Si^i)))}
\leq
 c(\|F\|_{L^p(\R_+; L^{q}_{\be_i}(\tilde\Om_i))}\ek
 \hspace{2cm}+\|\na U, P\|_{L^p(\R_+; L^{q}_0(\Om_0))}+\|U_t\|_{L^p(\R_+; (W^{1,q'}(\Om_0))')}),\; i=1,\ldots,m.
\end{array}
\end{equation}
 Thus, from \eq{4.18}-\eq{4.21} we get that
\begin{equation}
\label{E4.25}
\begin{array}{l}
\|u_0, u^{0}_t, \na^2 u^0, \na p^0\|_{L^p(\R_+; L^q(\Om_0))}\ek
\leq
 c(\|F,\na U, P\|_{L^p(\R_+; L^{q}(\Om_0))}
  +\|U_t\|_{L^p(\R_+; (W^{1,q'}(\Om_0))')}),\ek
 \|u^{i}_{t}, \na^2 u^i, \na p^i\|_{L^p(\R_+;  L^q_{\be_i}(\R; L^q(\Si^i)))}\leq
 c(\|F\|_{L^p(\R_+; L^{q}_{\be_i}(\tilde\Om_i))}\ek
 \hspace{2cm}+\|\na U, P\|_{L^p(\R_+; L^{q}(\Om_0))}+\|U_t\|_{L^p(\R_+; (W^{1,q'}(\Om_0))')}),\; i=1,\ldots,m.
\end{array}
\end{equation}
Note that $U=\sum_{i=0}^m u^i$, $P=\sum_{i=0}^m p^i$ in $\Om$. Therefore, by \eq{4.25} we have
\eq{4.33} and
\begin{equation}
\label{E4.32}
\begin{array}{l}
 \|U, U_{t}, \na^2 U, \na P\|_{L^p(\R_+;  L^q_{\tb}(\Om))}\leq
  c(\|F\|_{L^p(\R_+; L^{q}_{\tb}(\Om))}\ek
 \hspace{2cm}+\|\na U, P\|_{L^p(\R_+; L^{q}(\Om_0))}+\|U_t\|_{L^p(\R_+; (W^{1,q'}(\Om_0))')}).
\end{array}
\end{equation}

Consequently, it follows that the Stokes operator $A_{q,\tb}$ in $L^q_{\tb,\si}(\Om)$
 has maximal $L^p$-regularity for $1<p<\infty$ satisfying \eq{2.9}.

Thus, the proof of the Theorem \ref{T2.5} is complete, \qed
\par\bigskip\noindent
{\bf Acknowledgement:} Part of the work was done during the stay of the first author in the Institute
of Mathematics, AMSS, CAS, China under the support of
2012 CAS-TWAS Postdoctoral Fellowship, grant No. 3240267229.
This author is grateful to Prof. Ping Zhang for inviting me and to
CAS (Chinese Academy of Sciences) and TWAS (The World Academy
of Sciences) for providing me with the financial support.

%

\end{document}